\documentclass[11pt]{article}

\usepackage{latexsym}
\usepackage{amssymb}
\usepackage{color}
\usepackage{graphicx}

\newtheorem{Proposition}{Proposition}[part]

\newtheorem{Lemma}{Lemma}[part]

\def \R{\mathbb{R}}

\def \E{\mathbb{E}}

\def \P{\mathbb{P}}

\def \D{\mathbb{D}}

\def \ni{\noindent}

\title{\bf Multi-step Richardson-Romberg Extrapolation: \\ Remarks on Variance  Control and Complexity}

\author{{\sc Gilles Pag\`es} \thanks{Laboratoire de Probabilit\'es et Mod\`eles al\'eatoires, UMR~7599, Universit\'e Paris 6, case 188, 4,
pl. Jussieu, F-75252 Paris Cedex 5, France. E-mail:{\tt  gpa@ccr.jussieu.fr}}}

\date{January 17, 2007} 
\begin{document}

\maketitle

\begin{abstract}
We propose a multi-step Richardson-Romberg extrapolation method for the  computation
of expectations $\E f(X_{_T})$ of a diffusion $(X_t)_{t\in [0,T]}$ when the weak time
discretization error induced by the Euler scheme admits an expansion at an order $R\ge 2$.
The complexity of the estimator   grows as  
$R^2$ (instead of $2^R$ in the classical method) and its variance  is asymptotically controlled by considering
some consistent Brownian increments in the underlying Euler schemes. Some Monte Carlo simulations were carried with
path-dependent options (lookback, barrier) which support the conjecture that their weak time discretization error
also admits an expansion (in a different scale). Then an appropriate Richardson-Romberg extrapolation seems to
outperform the Euler scheme with Brownian bridge. 
\end{abstract}

\vspace{7mm}


\vspace{5mm}

\noindent {\bf MSC Classification (2000):} 65C05, 60H35, 65B99, 65C30.

\section{Introduction and preliminaries}

\setcounter{equation}{0}
\setcounter{Assumption}{0}
\setcounter{Theorem}{0}
\setcounter{Proposition}{0}
\setcounter{Corollary}{0}
\setcounter{Lemma}{0}
\setcounter{Definition}{0}
\setcounter{Remark}{0}

One considers a $d$-dimensional Brownian diffusion process $(X_t)_{t\in [0,T]}$ solution of the following
S.D.E.
\begin{equation}\label{EDS}
dX_t=b(t, X_t)dt+\sigma(t,X_t)dW_t,\qquad X_0=x.
\end{equation}
where $b:[0,T]\times \R^d\to \R^d$, $\sigma :[0,T]\times\R^d\to {\cal M}(d\times q)$ are continuous  
functions and
$(W_t)_{t\in [0,T]}$ denotes a $q$-dimensional Brownian motion defined on a filtered probability space
$(\Omega, {\cal A}, ({\cal F}_t)_{t\in [0,T]},  \P)$. We  assume that $b$ and $\sigma$ are
Lipschitz continuous in $x$ uniformly with respect to $t\!\in [0,T]$ and that $b(.,0)$ and $\sigma(.,0)$ are bounded over
$[0,T]$, see~\cite{KASC}.  In fact what we will simply need  is that, for every starting value $x\!\in \R^d$,    the Brownian Euler
scheme of~(\ref{EDS}) --
$i.e.$ the Euler scheme based on the increments of $W$ -- converges to   $X$ in every $L^p(\P)$,
$p\!\in(0,+\infty)$, where $X$ is the unique  strong
 solution  of the $SDE$  starting at
$x$. 
 
Except for some
very specific equations, it is impossible to process an exact simulation of the process
$X$ even at a fixed time
$T$ (by exact simulation, we mean writing
$X_{_T} =
\chi(U)$,
$U\sim U([0,1])$) (nevertheless,  when
$d=1$ and
$\sigma\equiv 1$, see~\cite{BEPARO},~\cite{BERO}). Consequently, to approximate
$\E(f(X_{_T}))$ by a Monte Carlo method, one needs to approximate $X$ by a process that can
be simulated (at least at a fixed number of instants). To this end one introduces the
stepwise constant Brownian Euler scheme $\bar X=(\bar X_{\frac{kT}{n}})_{0\le k\le n}$ with
step $\frac{T}{n}$  associated to the SDE. It is defined by
\[
\bar X_{t^n_{k+1}}= \bar X_{t^n_k}+b(t^n_k,\bar X_{t^n_k})\frac Tn+\sigma(t^n_k,\bar X_{t^n_k})
\sqrt{\frac Tn}\,U_{k+1},\qquad
\bar X_0=x,\quad k=0,\ldots,n-1,
\]
where $t^n_k= \frac{kT}{n},\, k=0,\ldots,n-1$ and $(U_k)_{1\le k\le n}$ denotes a sequence of i.i.d. ${\cal
N}(0;1)$-distributed random vectors given by
\[
U_k :=  \sqrt{\frac nT}\,(W_{t^n_k}-W_{t^n_{k-1}}),\quad
k=1,\ldots,n.
\]
Moreover, set for convenience $
\bar X_t := \bar X_{\underline t}$  where $\underline t=t^n_k$ if
$t\!\in[t^n_k,t^n_{k+1})$.

Then, it is classical background that
under the  regularity and growth assumptions on the coefficients $b$ and $\sigma$ mentioned
above,  $\sup_{t\in [0,T]}|X_t-\bar X_t|$ goes to zero in every $L^p(\P)$, $0<p<\infty$,  at a
$O(\frac{1}{\sqrt{n}})$-rate.

\smallskip However many authors in a long series of papers going back to the seminal papers by
Talay-Tubaro (\cite{TATU}) and Bally-Talay
(\cite{BATA1}, \cite{BATA2}), showed  under various assumptions on the
diffusion coefficients the existence of a vector space $V$  of   Borel 
functions $f:\R^d\to \R$ (bounded or with polynomial growth) for which one can expand the
``weak" time discretization error induced by the Euler scheme into a power series of
$\frac{1}{n}$. To be precise
\begin{equation}\label{Expansion}
({\cal E}^V_{_R})\hskip 0,75 cm \equiv \hskip 0,75 cm\forall\, f\!\in V,\qquad  \E(f(X_{_T}))= \E(f(\bar
X_{_T})) +
\sum_{k=1}^{R-1}
\frac{c_k}{n^k} + O(n^{-R})\hskip 0,75 cm 
\end{equation}
where the real constants $c_k$, $k= 1,\ldots, R-1$, do not depend on
the discretization parameter $n$ (see~\cite{TATU} for smooth functions $f$, 
see~\cite{BATA1} for bounded Borel functions  and uniformly hypo-elliptic
diffusions~(\footnote{In~\cite{BATA1}, the diffusion coefficients   are assumed to be ${\cal C}^\infty_b$  $i.e.$ infinitely  differentiable, bounded,  with bounded derivatives;  then $({\cal E}^V_{_R})$ holds for any $R$ if $\sigma$ is uniformly elliptic. As concerns $R=2$,    the
original  method of proof (see~\cite{TATU}) based on an approximation by the parabolic PDE involving the  infinitesimal
generator  works with finitely differentiable coefficients $b$, $\sigma$ and function $f$. Although not proved in full details, an  extension to Borel functions {\em with
polynomial growth} is mentioned in~\cite{BATA1}. See also~\cite{GOMU},~\cite{KOPE} for different approaches.})). In 
 ~\cite{BATA2}  is established the convergence of the p.d.f. of the Euler scheme at time $T$
toward that of $X_{_T}$. See also the recent work~\cite{GUY} for an extension to tempered
distributions. 
 Usually, this expansion is established in full
details  for $R=2$ and is known as the (standard) Richardson-Romberg extrapolation. However, up to additional
technicalities (and   smoothness assumptions on the coefficients)   the  expansion   holds true for larger values of
$R$ or even for every integer $R\ge2$. 

Furthermore, note that, the resulting vector space $V$ of ``admissible" functions is always  characterizing for
the
$\P$-$a.s.$ equality in the following sense: for every pair of $\R^d$-valued random vectors $X$, $Y$ defined on a
probability space
$(\Omega,{\cal A}, \P)$, 
\[
(\forall\, f\!\in V,\; f(X)=f(Y) \quad \P\mbox{-}a.s.)\Longrightarrow (X=Y\quad \P\mbox{-}a.s.)
\]
From now on, we always   assume that $V$ has this property. 

Sharp rates of weak convergence ($i.e.$ an expansion with $R=2$)
have been established for more general functionals
$F((X_t)_{t\in[0,T]})$ depending on the whole trajectory of the
diffusion. One uses the standard (stepwise constant) or the
continuous Euler scheme to approximate $\E F(X)$. These
contributions are usually motivated by the pricing of some classes
of path-dependent European options like~Asian options (~\cite{LATE}),
lookback options (see~\cite{SETO}~) or barrier options
(see~\cite{GOB}). These sharp rates are $c/n$ or $c/\sqrt{n}$ and one
can reasonably conjecture that these rates can be extended into some
expansions in an appropriate polynomial scale. The $R$-$R$
extrapolation can be in such situations a way to support some
conjectures concerning these expansions (and an heuristic way to
speed up the Monte Carlo simulations).

\medskip
In Section~\ref{RombStand}, we first briefly recall what the
original Richardson-Romberg ($R$-$R$) extrapolation is (when $R=2$).
Then, we prove that  the choice of consistent Brownian increments in
the two involved Euler schemes is optimal. In
Section~\ref{MultiRomb}, we propose a multi-step $R$-$R$
extrapolation at order $R$ and we show that choosing consistent
Brownian increments  still preserves the variance of the estimator
for  continuous functions or even Borel functions (under some
uniform ellipticity assumption). We point out that, however,
consistent increments is however not an optimal choice in general
for a given function $f$. Then, in Section~\ref{Complex} we analyze
the complexity of our approach and compare it to that of iterated
$R$-$R$ extrapolations.  In Section~\ref{Num}, we   provide for the
first values of interest ($R=2,3,4, 5$) some tables to simulate the
needed consistent Brownian increments at this order. In
Section~\ref{BSVanil}, we illustrate the efficiency of this
consistent $R$-$R$ extrapolation with $R=3,4$ by pricing vanilla
options in  a high volatility B-S model. In Section~\ref{BSpathdep},
we make some further tests on partial lookback and Up \& Out barrier
options which support the existence of expansions of the error.

\medskip
\noindent{\sc Notations:}
%$\bullet$
$A^*$ is for transpose of the matrix $A$. $|u|$ denotes the canonical Euclidean norm on $\R^q$, $q\ge 1$.
\section{The standard Richardson-Romberg extrapolation}\label{RombStand}

\setcounter{equation}{0}
\setcounter{Assumption}{0}
\setcounter{Theorem}{0}
\setcounter{Proposition}{0}
\setcounter{Corollary}{0}
\setcounter{Lemma}{0}
\setcounter{Definition}{0}
\setcounter{Remark}{0}
Assume that $({\cal E}^V_{_2})$ holds. Let $f\!\in V$ where $V$ denotes a vector space of continuous functions with
linear growth. The case of non continuous functions is investigated in the next section. For notational convenience
we set $W^{(1)}=W$ and $X^{(1)}:= X$. A regular Monte Carlo simulation based on 
$M$ independent copies
$(\bar X^{(1)}_{_T})^m$, $m=1,\ldots,M$, of the Euler scheme $\bar X^{(1)}_{_T}$ with step $T/n$ 
induces the following  global (squared) quadratic error
\begin{eqnarray}
\nonumber \|\E(f(X_{_T})) \!- \!\frac 1M\! \sum_{m=1}^Mf((\bar X^{(1)}_{_T})^m)\|^2_{_2}\!&\!=\!&\!|\E(f(X_{_T}))-
\E(f(\bar X^{(1)}_{_T}))|^2\\
\nonumber  &&+
\|\E(f(\bar X^{(1)}_{_T})) \!- \!\frac 1M\!
\sum_{m=1}^Mf((\bar X^{(1)}_{_T})^m)\|^2_{_2}\\
\label{L2error}\!&\!=\!&\!  \frac{c^2_1}{n^2}+ \frac{{\rm
Var}(f(\bar X^{(1)}_{_T}))}{M}  +O(n^{-3}).
\end{eqnarray}

This quadratic error bound~(\ref{L2error})  does not take fully advantage of the above expansion~$({\cal E}_{_2})$.
To take advantage of the expansion, one needs to make an $R$-$R$ extrapolation. In that framework (originally
introduced in~\cite{TATU}) one considers  a second Brownian Euler scheme, this time of the solution $X^{(2)}$  of a
``copy" of Equation~(\ref{EDS}) written with respect to a second Brownian motion $W^{(2)}$ defined on the same
probability space $(\Omega, {\cal A}, \P)$. In fact, one may chose this Brownian motion by enlarging $\Omega$ if necessary. This second Euler scheme has
{\em    a twice smaller}  step
$\frac{T}{2n}$ and is denoted
$\bar X ^{(2)}$. Then, assuming $({\cal E}^V_{_2})$to be more precise,  
\[
\E(f(X_{_T}))=\E(2f(\bar X^{(2)}_{_T})-f(\bar X^{(1)}_{_T}))  -\frac 12 \frac{c_2}{n^2} +O(n^{-3})
%+O(n^{-2}).
\]
Then, the new global (squared) quadratic error  becomes
\begin{equation}\label{L2error2}
\|\E(f(X_{_T}))- \frac 1M \sum_{m=1}^M 2f((\bar X ^{(2)}_{_T})^m)-f((\bar X_{_T})^m)\|_{_2}^2 = \frac{c^2_2}{4n^4}+ \frac{{\rm
Var}(2f(\bar X ^{(2)}_{_T})-f(\bar X_{_T} ^{(1)}))}{M} +O(n^{-5}).
\end{equation}

The structure of this quadratic error suggests the following question: is it possible to reduce the (asymptotic) time  discretization error
{\em  without increasing the Monte Carlo error? To what extend is it possible to control the variance term
${\rm Var}(2f(\bar X ^{(2)}_{_T})-f(\bar X^{(1)}_{_T}))$?}

\smallskip
%In order to answer this question in full generality we will assume that the original Euler scheme
%with step $\frac Tn$ is related to a (first) Brownian motion $W= W^{(1)}$ $i.e.$ $U_k =
%U^{(1)}_k:=\sqrt{\frac Tn}(W^{(1)}_{t^n_k}-W^{(1)}_{t^n_{k-1}})$ whereas the second scheme (with step
%$\frac{T}{2n}$) is related to a second Brownian $W^{(2)}$ (defined on the same probability space as $W^{(1)}$)
%$i.e.$
%\begin{equation}\label{Eulerhalf}
%\bar X ^{(2)}_{t^{2n}_{k+1}}= \bar X ^{(2)}_{t^{2n}_k}+b(t^{2n}_k,\bar
%X^{(2)}_{t^{2n}_k})\frac{T}{2n}+\sigma(t^{2n}_k,\bar X ^{(2)}_{t^{2n}_k})
%\sqrt{\frac{T}{2n}}\,U^{(2)}_{k+1},\;
%\bar X ^{(2)}_0=x,\; k=0,\ldots,2n-1,
%\end{equation}
%with $U^{(2)}_k:=\sqrt{\frac{2n}{T}}(W^{(2)}_{t^{2n}_k}-W^{(2)}_{t^{2n}_{k-1}})$.

\smallskip
It is clear that $\sup_{t\in[0,T]}|\bar X ^{(i)}_{t}-X^{(i)}_{t}|$ converges to $0$ in every $L^p(\P)$,
$0<p<\infty$, $i=1,2$. 
%Once
%extended over $[0,2T]$ the Euler scheme $\bar X$ converges similarly toward the solution $(X^{(1)}_t)_{t\in [0,2T]}$
%of~(\ref{EDS}) on $[0,2T]$ (with $W= W^{(1)}$). 
Consequently
\[
\sup_{t\in[0,T]}|(\bar X^{(1)}_t,\bar X ^{(2)}_{t})-(X^{(1)}_t,X^{(2)}_{t})|\stackrel{L^p(\P) \, \& \,
a.s.}{\longrightarrow}0\quad\mbox{ as }\quad n\to \infty.
\]
In particular (keep in mind $f$ has at most polynomial growth),
\[
 {\rm Var}(2f(\bar X ^{(2)}_{_T})-f(\bar X^{(1)}_{_T})) \longrightarrow {\rm
Var}(2f(X^{(2)}_{_{T}})-f(X^{(1)}_{_T}))\qquad
\mbox{ as }\qquad n\to \infty.
\]
 
Then, straightforward computations show that

\begin{equation}\label{Var2}
{\rm Var}\left(2f(X^{(2)}_{_{T}})-f(X^{(1)}_{_T})\right)= 5 \,\E(f(X^{(1)}_{_{T}}))^2 - 4
\,\E f(X^{(1)}_{_T})f(X^{(2)}_{_{T}}) -\left(\E(f(X^{(1)}_{_{T}}))\right)^2
\end{equation}
where
\begin{eqnarray} \label{Cov} 
\E \left(f(X^{(1)}_{_T})f(X^{(2)}_{_{T}})\right)&\le &\|f(X^{(1)}_{_T})\|_{_2}
\|f(X^{(2)}_{_{T}})\|_{_2}= \|f(X^{(1)}_{_T})\|_{_2}^2\hskip 0,5  cm
\end{eqnarray}
by  Schwarz's Inequality since $X^{(1)}_{_T}\stackrel{d}{=} X^{(2)}_{_T}$. 

Consequently, minimizing the variance term amounts to maximizing $\E( f(X^{(1)}_{_T})f(X^{(2)}_{_{T}}))$. It
follows from the equality case in Schwarz's Inequality that  the  equality in the left hand side of~(\ref{Cov}) holds  iff
$f(X^{(2)}_{_T})= \lambda f(X^{(1)}_{_T})$ $\P$-$a.s.$ for some $\lambda\!\in \R_+$ or $f(X^{(1)}_{_T})=0$ $\P$-$a.s.$. But since $X^{(1)}_{_T}$ and
$X^{(2)}_{_T}$ have the same distribution, this   implies $f(X^{(2)}_{_T})= f(X^{(1)}_{_T})$ $\P$-$a.s.$.

Finally minimizing the asymptotic variance term for every $f\!\in V$ is possible if and only if 
\begin{equation}\label{X1=X2}
X^{(2)}_{_T}= X^{(1)}_{_T}\quad \P\mbox{-}a.s.
\end{equation}

A sufficient (and sometimes necessary, see Annex) condition to ensure~(\ref{X1=X2}) clearly is   that
\[
 W^{(2)} = W^{(1)}.
\]
This means that the white noise $(U^{(2)}_k)_{1\le k\le n}$ of the Euler scheme $\bar X^{(2)}$ satisfies 
\[
U^{(2)}_k = \sqrt{\frac{2n}{T}}(W^{(1)}_{\frac{kT}{2n}}-W^{(1)}_{\frac{(k-1)T}{2n}}),\quad k=1,\ldots,2n,
\]
so that, from a simulation point of view, one needs to simulate $2n$ i.i.d. copies  $U^{(2)}_k$ of
the normal distribution and then sets
\[
U^{(1)}_k = \frac{U^{(2)}_{2k}+U^{(2)}_{2k-1}}{\sqrt{2}},\; k=1,\ldots, n.
\]

Note that all what precedes (as well as all what follows in fact) can be
extended to continuous functionals $F(X)$ for which   an expansion
similar to~(\ref{Expansion}) holds (see $e.g.$ the pricing of
  Asian options in~\cite{LATE} or Section~\ref{Num}) or  for different dynamics like SDE with delay studied
in~\cite{CLKOLA}.

\smallskip
This result provides an  (at least asymptotic) positive answer  to a question raised in~\cite{GLA} (see p. 361): Self-consistent
Brownian increments do minimize the (asymptotic) variance in the
$R$-$R$ extrapolation. In the next section we show that it is possible to design some multi-step $R$-$R$
extrapolations with a reasonable complexity for which the control of the variance is still preserved.   

As concerns variance control, note that if one proceeds using two {\em independent} sequences of Gaussian
white  noises $U^{(1)}$ and $U^{(2)}$, the expansion~(\ref{Var2}) yields  
\[
 {\rm Var}\left(2f(X^{(2)}_{_{T}})-f(X^{(1)}_{_T})\right) = 5\, {\rm Var}\left(f(X^{(1)}_{_T})\right)
\]
Such a choice is then the worst possible one. It induces an increase of the Monte Carlo simulation by a factor $5$.
\section{Multi-step Richardson-Romberg extrapolation}\label{MultiRomb}

\setcounter{equation}{0}
\setcounter{Assumption}{0}
\setcounter{Theorem}{0}
\setcounter{Proposition}{0}
\setcounter{Corollary}{0}
\setcounter{Lemma}{0}
\setcounter{Definition}{0}
\setcounter{Remark}{0}

Iterated $R$-$R$ extrapolations are usually not
implemented, essentially because their numerical efficiency (when the white noises are independent) becomes less and
less obvious. The first reason is the increase of the  complexity, the second 
one is the ``explosion" of the variance term (especially when implemented with independent Brownian motions) and the third one is
the absence of control of the coefficients
$c_k$ as $k$ increases.
 In what follows we   propose a solution to the first two problems: we propose a multi-step $R$-$R$
extrapolation with consistent Brownian increments. Doing so, we will control the variance in the Monte
Carlo error term and limit the increase of the complexity of the procedure.  
\smallskip
One proceeds as follows: let  $R\ge 2$ be an integer.  We wish to obtain a time
discretization error behaving like $n^{-R}$ as $n\to \infty$. To this end, we introduce
$R$ Gaussian Euler schemes $\bar X^{(r)}= (\bar X^{(r)}_{\frac{kT}{rn}})_{0\le k\le r\,n} $ with step
$\frac{T}{r\,n}$, $r=1,\ldots,R$. Each Euler scheme is designed using a Gaussian white noise
$(U^{(r)}_k)_{1\le k\le r\, n}$ obtained from the increments of a standard Brownian motion
$W^{(r)}$ by setting
\[
U^{(r)}_k:= \sqrt{\frac{r\,n}{T}}\left(W^{(r)}_{\frac{kT}{rn}}-W^{(r)}_{\frac{(k-1)T}{rn}}\right),\; k=1,\ldots, r\,n.
\]
The Brownian motions $W^{(1)},\ldots, W^{(R)}$ are all defined on the same probability space.
We will come back further on about the practical simulation of these increments. Our
``meta-assumption" $({\cal E}^V_{_{R+1}})$ implies that for a  function
$f\!\in V$, one has
\[
 \E(f(X_{_T})) = \E(f(\bar X^{(r)}_{_T})) + \sum_{\ell=1}^{R-1} \frac{c_\ell}{r^\ell}\,\frac{1}{n^\ell} +
\frac{c_{_R}}{r^R}\frac{1}{n^{R}}(1+O(1/(rn))),\qquad
r=1,\ldots,R.
\]
One defines the $R\times (R-1)$ matrix
\[
A=\left[\frac{1}{r^{\ell}}\right]_{1\le r\le R, 1\le \ell \le R-1}
%= {\rm Vandermonde} \left(\frac 1r, r=1,\ldots,
%R\right).
\]
and $\mbox{\bf 1}$ as the unit (column) vector of $\R^{R}$. Then, the above system of equations reads
\begin{equation}\label{Systeme}
\E(f(X_{_T})) \mbox{ \bf 1} = \left[\begin{array}{c} \vdots\\ \E(f(\bar X^{(r)}_{_T}))\\
\vdots\end{array}\right]_{1 \le r\le R}\hskip -0,5cm+A \left[\begin{array}{c} \vdots\\ \frac{c_r}{n^r}\\
\vdots\end{array}\right]_{1 \le r\le R-1} \hskip -0.9cm+ c_{_R}\left[\begin{array}{c}
\vdots\\\frac{1}{n^Rr^{R}}(1+O(\frac{1}{r\,n}))
\\
\vdots\end{array}\right]_{1 \le r\le R}.
\end{equation}
 Now, let $\alpha \!\in \R^R$ satisfying
$$
\alpha^*\,\mbox{\bf 1}=1\qquad \mbox{ and } \qquad \alpha^*\, A =0.
$$
This reads 
\begin{equation}\label{equalpha}
\alpha^* \, \widetilde A = e_1
\end{equation}
where $e_1= (1,0,\ldots,0)^*$ is the first (column) vector of the canonical basis of $\R^R$ and
\[
\widetilde A = \left[\frac{1}{r^{\ell-1}}\right]_{1\le r\le R, 1\le \ell \le R}
= {\rm Vandermonde} \left(\frac 1r, r=1,\ldots,
R\right).
\]
Multiplying~(\ref{Systeme}) on the left by $\alpha^*$ yields
\begin{eqnarray}
\label{RmultiRomb}\E(f(X_{_T})) &=& \E\left(\sum_{r=1}^R \alpha_r f(\bar X^{(r)}_{_T})
\right) +\frac{\widetilde c_{_R}}{n^R}\left(1+ O(1/n)\right)\hskip 3cm\\
\nonumber\mbox{where  }\hskip 3 cm  \widetilde c_{_R} &= & c_{_R}
\sum_{r=1}^R\frac{\alpha_r}{r^{R}}.
\end{eqnarray}

Solving the Cramer linear system~(\ref{equalpha}) yields after some elementary computations based on the Vandermonde determinant
 the following lemma.

\begin{Lemma} For every integer $R\ge 1$,
\begin{equation}\label{alpha}
 \alpha_r =\alpha(r,R):=  (-1)^{R-r} \frac{r^R}{r!(R-r)!},\quad
1\le r\le R.
\end{equation}
\end{Lemma}

\noindent {\bf Remarks.} $\bullet$  As a consequence of~(\ref{alpha}), one shows that for every integer $R\ge 1$,
\[
\sum_{r=1}^R\frac{\alpha_r}{r^{R}}= \sum_{r=0}^{R-1} \frac{(-1)^r}{r!} \frac{1}{(R-r)!}= \frac{(-1)^{R-1}}{R!}
\]
(this follows from the comparison of the coefficient of $x^R$ in the expansion $e^{x}e^{-x}=1$) so that
$$
|\widetilde c_{_R}|= \frac{|c_{_R}|}{R!}.
$$
This bound emphasizes that this multi-step  extrapolation {\em does not}  ``hide"   an ``exploding" behaviour   in the
first non vanishing  order of the expansion. In fact it has a damping effect on this coefficient $c_{_R}$. 
However, the fact that we have in practice no control on the coefficient $c_{_R}$ (and its successors) at a reasonable cost  remains and is not  overcome by this approach. 

\smallskip \noindent $\bullet$ If ${\cal E}^V_{_R}$ holds, then the above computations show that~(\ref{RmultiRomb}) becomes
\[
\E(f(X_{_T})) = \E\left(\sum_{r=1}^R \alpha_r f(\bar X^{(r)}_{_T})
\right) +O(n^{-R}).
\]
 
\smallskip
 Now we pass to the choice of the Brownian motions
$W^{(r)}$ in order to preserve  the variance of such a combination. 

\subsection{The case of continuous functions with polynomial growth}\label{SecCont}
In this section, we assume that the  vector space $V$ is made  of   continuous functions  with
polynomial growth and is stable by product.  First note that for every $p>0$ 
\begin{equation}\label{ConvEulerR}
\sup_{t\in[0,T]}|(\bar X^{(1)}_{t},\ldots,\bar X^{(r)}_{t},\ldots, \bar
X^{(R)}_{t})-(X^{(1)}_t,\ldots,X^{(r)}_{t},\ldots, X^{(R)}_{t})|\stackrel{L^p(\P)\,\&\, a.s.}{\longrightarrow}0.
\end{equation}

Let $f\!\in V$. Then 
%(\footnote{In many situations, this convergence holds at a $O(1/n)$-rate
%since one can reasonably expect that the SDE associated to the process $(X^{(1)},\dots, X^{(r)},\dots, X^{(R)})$ satisfies $({\cal
%E}_{_R})$.})
\begin{equation}\label{VarConv}
{\rm Var} \left(\sum_{r=1}^R \alpha_r f(\bar X^{(r)}_{_T}) \right)\longrightarrow {\rm Var}\left(\sum_{r=1}^R
\alpha_r f( X^{(r)}_{_{T}})\right)\qquad\mbox{ as }\quad n\to \infty
\end{equation}

\noindent since  $f^2$ is  continuous with polynomial growth.  When 
\[
  W^{(r)} = W, \quad r=1,\ldots, R,
\]
then, $f( X^{(r)}_{_{T}})= f(X_{_{T}})$ for every $r\!\in \{1,\ldots,R\}$ so that 
\[
 {\rm Var}\left(\sum_{r=1}^R
\alpha_r f( X^{(r)}_{_{T}})\right) = \left(\sum_{r=1}^R \alpha_r\right)^2 {\rm Var}(f( X_{_{T}}))= {\rm
Var}(f( X_{_{T}})).
\]
This obvious remark shows that this  choice for the Brownian motions $W^{(r)}$ lead to a control of the
variance of the multi-step $R$-$R$ estimator by that of $f(X^{(1)}_{_T})$. However, the  
optimality of this choice turns out to be a less straightforward question. 

On the opposite, if one considers some mutually independent Brownian motions $W^{(r)}$ (which is in some way the
``laziest" choice from the programming viewpoint), then
\[
 {\rm Var}\left(\sum_{r=1}^R
\alpha_r f( X^{(r)}_{_{T}})\right) =  \left(\sum_{r=1}^R \alpha^2_r\right) {\rm Var}(f( X_{_{T}})) 
\] 
 having in mind that
\[
\sum_{r=1}^R \alpha^2_r\ge \alpha^2_{_R} = \left(\frac{R^R}{R!}\right)^2\sim\frac{e^{2R}}{2\pi R}\quad \mbox{ as }\quad R\to
\infty.
\]
The practical aspects of the above  analysis can be  summed up   as follows.
\begin{Proposition}\label{PropCont} Let $T>0$. One considers the SDE~(\ref{EDS}). Assume that
 on every filtered probability space $(\Omega,{\cal
A},({\cal F}_t)_{t\in [0,T]}, \P)$ on which exists a standard ${\cal F}_t$-Brownian motion $W$,  
~(\ref{EDS})  has a   strong solution   on
$[0,T]$  and that its Euler scheme with step $T/n$ (with Brownian increments) converges in  every $L^p(\P)$ for the sup-norm over $[0,T]$ for every
$p\!\in[1,\infty)$.
% Assume the same holds true for
%the $SDE$'s  with coefficients $\lambda \, b$ and $\sqrt{\lambda} \, \sigma$ for every $\lambda >0$. 
Furthermore assume
that~(\ref{EDS})  admits an expansion~$({\cal E}^V_{_R})$  at an
order
$R\ge 1$ for a characterizing vector space  $V$ of   continuous functions $f$ with polynomial growth, stable by
product.

Let   $\alpha\!\in \R^R$ be defined by~(\ref{alpha}). For every
$r\!\in\{1,\ldots,R\}$ set
\begin{equation}\label{ConsistIncr}
U^{(r)}_k:= \sqrt{\frac{rn}{T}}\left(W_{\frac{kT}{rn}}-W_{\frac{(k-1)T}{rn}}\right),\;
k\ge 1.
\end{equation}
If $\bar X^{(r)}$ denotes an Euler scheme of~(\ref{EDS}) with step $\frac{T}{rn}$ associated to the
Gaussian white noise $(U^{(r)}_k)_{k\ge 1}$, then
\[
\E\left(\!\sum_{r=1}^R\!\alpha_r f(\bar X^{(r)}_{_T})\!\right) \!= \!\E(f(X_{_T})) +O(n^{-R})\;\mbox{ and }\;
\lim_n {\rm Var}\left(\!\sum_{r=1}^R\!\alpha_r f(\bar X^{(r)}_{_T})\!\right)= {\rm Var}(f(X_{_T})).
\]
\end{Proposition}

\noindent {\bf Remark.}   A first extension of this result to non continuous Borel functions can be made easily
as follows: Convergence~(\ref{ConvEulerR}) implies that $(\bar X^{(1)}_{_T},\ldots,\bar X^{(r)}_{_T},\ldots, \bar
X^{(R)}_{_T})$ weakly converges in $(\R^d)^R$  to $(X^{(1)}_{_T}, \ldots,X^{(1)}_{_T})$. Then, if $f$ is $\P_{X_{_T}}(dx)$-$a.s.$
continuous with polynomial growth, so is $(\sum_r \alpha_rf(x_r))^2 $ with respect to $\P_{(X^{(1)}_{_T}, \ldots,X^{(1)}_{_T})}$.
This implies that the  convergence of the  variances~(\ref{VarConv}) still holds true. 

%$\bullet$ The extension of the above result to possibly non continuous Borel functions would require to prove that
%$\E(f(\bar X^{(r)}_{_{T}})f(\bar X^{(r')}_{_{T}})) \to \E(f( X^{(r)}_{_{T}})f(X^{(r')}_{_{T}}))$ for every $r\neq r'$ which does not 
%straightforwardly  follow from the existing literature on the expansion of the weak error.
\subsection{The case of bounded Borel functions (elliptic diffusions)}
We no longer assume that $V$ is a subspace of continuous functions. As a counterpart, we make some more stringent assumptions on the
 coefficients $b$ and $\sigma$ of the diffusion~(\ref{EDS}), namely that they satisfy the following assumption:
\[
(UE) \equiv \left\{\begin{array}{ll}
(i)&b_i,\,\sigma_{ij} \!\in{\cal C}^\infty_b(\R^d,\R),\; i\!\in\{1,\ldots,d\},\;
\,j\!\in\{1,\ldots,q\}\\ *[.4em] (ii)&\exists\, \varepsilon_0>0 \hbox{ such that }\forall\, x\!\in \R^d,\; \sigma\sigma^*(x) \ge \varepsilon_0 I_d.
\end{array}\right.
\] 
(In particular this implies $q\ge d$). Then we know from~\cite{BATA1} (see also~\cite{GUY}) that $({\cal E}_{_R}^V)$ holds at any
order  
$R\ge 1$ with  
$V={\cal B}_b(\R^d,\R)$ where  ${\cal B}_b(\R^d,\R)$ denotes the set of bounded Borel functions $f:\R^d\to \R$. In that setting, the
following result holds
\begin{Proposition}\label{PropBor}
 If Assumption~$(UE)$ holds, then Proposition~\ref{PropCont} still holds at any order $R\ge 1$ with $V={\cal B}_b(\R^d,\R)$.
\end{Proposition}

\noindent {\bf Proof.} The conclusion of Proposition~\ref{PropCont} in Section~\ref{SecCont}  
remains valid for a subspace $V$ of possibly non continuous functions provided the convergence of
the variance term in Equation~(\ref{VarConv}) still holds for every $f\!\in V$. Developing the variance term as
follows 
\[
{\rm Var}\left(\sum_{r= 1}^R \alpha_r f(\bar X^{(r)}_{_T})\right)= \sum_{r,r'=1}^R \alpha_r\alpha_{r'} \E\left(f(
\bar X^{(r)}_{_T})f(\bar  X^{(r')}_{_T})\right)-\left(\sum_{r=1}^R \alpha_r
\E f(\bar X^{(r)}_{_T})\right)^2
\]
shows that it amounts to proving that
\begin{equation}\label{CovBor}
\forall\, r,\, r' \!\in \{1,\ldots,R\},\quad \E\left(f(\bar X^{(r)}_{_T})f(\bar X^{(r')}_{_T})\right) \longrightarrow \E\left(f(
X^{(r)}_{_T})f( X^{(r')}_{_T})\right) \; \mbox{ as } \; n\to \infty.
\end{equation}
In order to prove this convergence, we rely on the following lemma established in~\cite{BATA2}.
\begin{Lemma} If $(UE)$ holds then the distributions of  $X_{_T}$ and its Brownian Euler schemes
with step $T/n$ , starting at
$x\!\in \R^d$ are absolutely continuous with distributions $p_{_T}(x,y)\lambda_d(dy)$ and
$\bar p^n_{_T}(x,y)\lambda_d(dy)$ respectively. Furthermore, they satisfy  for every
$x,y\!\in \R^d$ and every $n\ge 1$, 
\[
\bar p^n_{_T}(x,y) = p_{_T}(x,y) + \frac{\pi_n(x,y)}{n}\quad \hbox{ with }\quad |\pi_n(x,y)| \le C_1e^{-C_2|x-y|^2}
\]
for a real constants $C_1,\,C_2>0$.
\end{Lemma}

Let $\mu$ denote the finite positive measure defined by $\mu(dy)= (p_{_T}(x,y) + C_1e^{-C_2|x-y|^2})\lambda_d(dy)$. The set ${\rm Lip}_b(\R^d,\R)$ of
bounded Lipschitz functions is everywhere dense in $L^1_{\R}(\R^d,\mu)$. Consequently, there exists a sequence $(f_k)_{k\ge 1}$ of functions of
${\rm Lip}_b(\R^d,\R)$ such that $\lim_k \|f-f_k\|_{_1}=0$. Furthermore, since one can replace each $f_k$ by $((-\|f\|_{_\infty})\vee f_k)\wedge
\|f\|_{_\infty}$, one can assume without loss of generality that $\|f_k\|_{_\infty}\le \|f\|_{_\infty}$.   Then one gets
\begin{eqnarray*}
\left|\E\left(f(
\bar X^{(r)}_{_T})f(\bar  X^{(r')}_{_T})- f(X^{(r)}_{_T})f(X^{(r')}_{_T})\right)\right|&\le& \left|\E\left(f(
\bar X^{(r)}_{_T})f(\bar  X^{(r')}_{_T})-f_k(\bar
X^{(r)}_{_T})f_k(\bar X^{(r')}_{_T})\right)\right|\\
&&+ \left|\E\left(f_k(
\bar X^{(r)}_{_T})f_k(\bar  X^{(r')}_{_T})- f_k(
X^{(r)}_{_T})f_k(X^{(r')}_{_T})\right)\right|\\
&&+\left|\E\left(f_k(
X^{(r)}_{_T})f_k(X^{(r')}_{_T})-f(X^{(r)}_{_T})f(X^{(r')}_{_T})\right)\right|\\
&\le&  2\|f\|_{_\infty}\max_{1\le r\le R} \E\left(|f-f_k|(\bar X^{(r)}_{_T})\right)\\
&&+ 2\|f\|_{_\infty}[f_k]_{\rm Lip}\max_{1\le r\le R}\E| X^{(r)}_{_T}-\bar X^{(r)}_{_T}|\\
&&+2\|f\|_{_\infty}\max_r \E\left(|f-f_k|(X^{(r)}_{_T})\right).
\end{eqnarray*}

Now, for every $k\ge 1$ and every $r\!\in\{1,\ldots, R\}$, 
\[
\E\left(|f-f_k|(\bar X^{(r)}_{_T})\right)+\E\left(|f-f_k|(X^{(r)}_{_T})\right) \le 2\|f-f_k\|_{L^1(\mu)}
\]
so that, for  every $r,\,r'\!\in\{1,\ldots, R\}$ and every  $k\ge 1$,
\[
\limsup_n \left|\E\left(f(
\bar X^{(r)}_{_T})f(\bar  X^{(r')}_{_T})- f(X^{(r)}_{_T})f(X^{(r')}_{_T})\right)\right|\le  2\|f\|_{_\infty}\|f-f_k\|_{L^1(\mu)}
\]
which in turn implies~(\ref{CovBor}) by letting $k$ go to infinity.$\qquad_{\diamondsuit}$ 

\subsection{About the optimality of the Brownian specification when $R\ge 3$}

We proved in Section~\ref{RombStand} that  for the standard $R$-$R$ extrapolation ($R=2$), the choice of the same
underlying Brownian motion $W$ for both schemes is optimal $i.e.$ leads to the lowest possible asymptotic variance
for any admissible function $f$. In fact this is the only case of optimality: the choice of consistent Brownian
increments in the  higher order
$R$-$R$ extrapolation  is never optimal in terms of variance reduction when $R\ge 3$.

 Let
$f\!\in V$ denote a (fixed) function in $V$ such that ${\rm Var}(f(X^{(1)}_{_T}))>0$.  One may assume without
loss of generality that  
$$
{\rm Var}(f(X^{(1)}_{_T}))=1.
$$ 
Let $S^f := [{\rm
Cov}(f(X^{(r)}_{_T}),f(X^{(r')}_{_T}))]_{1\le r,r'\le R}$ denote the covariance matrix of the variables $f(X^{(r)}_{_T})$. Since  ${\rm
Var}(f(X^{(r)}_{_T}))=1$, $r=1,\ldots,R$, the minimization of the variance   of the estimator 
$\alpha^*\,[f(X^{(r)}_{_T})]_{1\le r\le R}$ is ``lower-bounded" in a natural way by  
  the following problem   
\begin{equation}\label{minS}
\min \left\{\alpha^* S\alpha, \; S\!\in {\cal S}^+(d,\R), \,S_{ii}=1 \right\}.
\end{equation} 
Choosing like in Propositions~\ref{PropCont} and~\ref{PropBor},   $f(X^{(r)}_{_T})=f(X_{_T})$ $\P$-$a.s.$,
$r=1,\ldots,R$    corresponds to
$S^f=S  = 
\mbox{\bf 1}\,\mbox{\bf 1}^*$ ($i.e.$ $S^f_{ij}=1$, $1\le i,j\le R$). 

%Using Cholevsky's decomposition, $S^f =T^*T$, $T$ upper triangular,  shows that this is equivalent to 
%\begin{equation}\label{minT}
%\min \left\{|T\alpha|^2,\; T=[t_{ij}]  \hbox{ lower triangular}, \sum_{i=1}^j t^2_{ij}=1,\, j=1,\ldots,R \right\}.
%\end{equation}

We use the term  ``lower-bounding" to emphasize that for a given $f\!\in V$, there are possibly admissible matrices
in~(\ref{minS}) which cannot be covariance matrices for $(f(X^{(r)}_{_T}))_{1\le r\le R}$.

\begin{Proposition} $(a)$ When $R\ge 3$, 
\[
\min \left\{\alpha^* S\alpha, \; S\!\in {\cal S}^+(R,\R), \,S_{ii}=1 \right\}=0,
\]
hence, the unit matrix $S_1:=\mbox{\bf 1}\,\mbox{\bf 1}^*$ is never a solution to the minimization
problem~(\ref{minS}).

\smallskip
\noindent  $(b)$ When $R=3$, if $f(X_{_T})= \E f(X_{_T}) + {\rm StD}(f(X_{_T}))\varepsilon$,
$\P(\varepsilon=\pm 1)=1/2$,  then the extrapolation formula stands as an exact quadrature formula. 
\end{Proposition}

\noindent {\bf Remark.} Although it has no practical interest for applications, the situation described in $(b)$
may happen: assume the drift $b$ is odd, the diffusion coefficient $\sigma$ is even and $f$ is the sign function
$f(x) = {\rm sign}(x)$. When $R\ge 4$, for a given function $f$,    the solution(s) of
the abstract minimization problem~(\ref{minS})   do not correspond in general to   an
admissible covariance matrix for
$(f(X^{(r)}_{_T}))_{1\le  r\le R}$ (see below the proof for a short discussion on simple examples).   Furthermore, even when it
happens to be the case, the variance control provided by
$\mbox{\bf 1}\,\mbox{\bf 1}^*$ is obviously more straightforward  than a
numerical search of   the appropriate underlying covariance structure of the Brownian motions
$(W^{(r)})_{1\le  r\le R}$ driving the diffusions. However, one must keep in mind that multi-step $R$-$R$ extrapolation leaves some  degrees of
freedom to some variance reduction method as soon as $R\ge3 $ and the opportunity in terms of computational cost to design an
online optimization procedure (depending on the function $f$) may be interesting in some cases.  The appropriate approach is then
an online recursive stochastic approximation method somewhat similar to   that introduced in~\cite{AR} for online variance
reduction. This is beyond the scope of the present paper.

%It turns out the lowest achievable variance $\alpha S^f\alpha \!\in (0,1)$ and that\dots

\bigskip
\noindent {\bf Proof.} $(a)$ Assume $R\ge 4$. We will show that the minimum in~(\ref{minS}) is in fact $0$. Let $I^{\pm}:=
\{i\!\in\{1,\ldots,R\}\,|\,
\alpha_i\begin{array}{c} >\\ <
\end{array}0\}$. Let
$\rho\!\in (0,1)$ and let $A(\rho):=\sqrt{(1-\rho) I_R+\rho
\mbox{\bf 1}\mbox{\bf 1}^*}$ denote the  symmetric square root of the positive definite matrix  $(1-\rho) I_R+\rho
\mbox{\bf 1}\mbox{\bf 1}^*$. Then set for every $i\!\in I^+$,  $C_i = A(\rho) e_i$ where $(e_1,\ldots, e_R)$
still denotes the canonical basis of
$\R^R$ and 
\[
\theta_i := \frac{\alpha_i}{\sum_{j\in I^-}|\alpha_j|},\qquad i\!\in I^+.
\] 
First note that $|C_i|^2 = \sum_{j}(A(\rho)_{ji})^2 = (A(\rho)^2)_{ii}=1$ and $(C_i|C_j)=\rho$, $i\neq j$. Then, 
\[
\sum_{i\in I^+} \theta_i = \frac{\sum_{i\in I^+}|\alpha_i|}{\sum_{j\in I^-}|\alpha_j|}= \frac{1+\sum_{j\in I^-}|\alpha_j|}{\sum_{j\in
I^-}|\alpha_j|}>1. 
\]
On the other hand
\[
\sum_{i\in I^+} \theta^2_i = \frac{\sum_{i\in I^+}\alpha_i^2}{(\sum_{j\in I^-}|\alpha_j|)^2}=\frac{\sum_{i\in I^+}\alpha_i^2}{(\sum_{i\in
I^+}\alpha_i-1)^2} = \frac{\sum_{i\in I^+}\alpha_i^2}{\sum_{i\in I^+}\alpha_i^2 +1+\sum_{i\in I^+}\alpha_i(\sum_{j\neq i,j\in I^+}\alpha_j-2)} <1
\]
since for every $i\!\in I^+$, $\sum_{j\neq i,j\in I^+}\alpha_j\ge \min(\alpha_{R-2}, \alpha_R)>2$ provided $R\ge 4$. Consequently
the  function $g:[0,1]\to \R$ defined by $g(\rho):=|\sum_{i\in  I^+}\theta_i C_i|^2$ satisfies $g(0)=\sum_{i\in I^+} \theta^2_i<1$ and
$g(1) = (\sum_{i\in I^+} \theta_i)^2>1$ so that there exists a real number $\rho_0\!\in (0,1)$ such that $g(\rho_0)=1$. From now on assume that
$\rho=\rho_0$. 

Set $C_i= C:=\sum_{j\in I^+} \theta_j C_j$, for every $i\!\in I^-$  so that $|C|=1$  and 
\[
\sum_{r=1}^R \alpha_i C_i=(\sum_{i\in I^-} \alpha_i)C + \sum_{i\in I^+}\alpha_i C_i=0.
\]
Consequently the symmetric nonnegative matrix 
\begin{equation}\label{Smini}
S:= [C_1\cdots C_{_R}]^*\,[C_1\cdots C_{_R}]
\end{equation}
satisfies $S_{ii}=1$ for every $i\!\in\{1,\ldots,R\}$ and $S \alpha=0$.

 When $R=3$, one checks that $S^f =
u\,u^*$,  $u^*= [1 \,-1\,-1]$,  satisfies $S^f \alpha=0$ since $\alpha_2+\alpha_3 =1/2= \alpha_1$.  

\medskip 
\noindent  $(b)$  The situation for $R=3$ 
corresponds to    $f(X^{(2)}_{_T})= f(X^{(3)}_{_T})=2\,\E(f(X^{(1)}_{_T}))-f(X^{(1)}_{_T})$ so
that
$f(X^{(1)}_{_T})+f(X^{(2)}_{_T})= 2\,\E f(X^{(1)}_{_T})$ $\P$-$a.s.$. The same identity holds
for $f^2\!\in V$, so that $f(X^{(1)}_{_T})\times f(X^{(2)}_{_T})= \E f(X^{(1)}_{_T})f(X^{(2)}_{_T})$ $\P$-$a.s.$
(use $ab= \frac 12((a+b)^2-a^2-b^2)$). Consequently $f(X^{(1)}_{_T})$ and $f(X^{(2)}_{_T})$ are the zeros of
$z^2-uz+v$ for some deterministic real constants $u$ and $v$. The conclusion follows once noticed that
$f(X^{(1)}_{_T})\stackrel{d}{=} f(X^{(2)}_{_T})$. 
$\quad_{\diamondsuit}$

\bigskip
\noindent {\sc Two toy situations:} $(a)$ Assume -- which is of no numerical interest -- that $X^{(r)}_{_T} =
W^{(r)}_{_T}$,
$r=1,\ldots,R$, where the covariance of the Brownian motions $W^{(r)}$ is given by~(\ref{Smini}). Then $S\alpha
=0$ so that the multi-step $R$-$R$ extrapolation is an exact quadrature formula (to compute $0$\dots). 

\medskip
\noindent $(b)$ If one considers some correlated geometrical Brownian motions $X_t^{(r)}= e^{-\frac{\sigma^2t}{2} +\sigma
W^{(r)}_t}$ $r=1,\ldots,R$, where $(W^{(1)},\ldots, W^{(R)})$ has $C_{_W}=[c_{rr'}]_{1\le r,r'\le R}$ as a covariance matrix (at
time
$t=1$), then the covariance matrix of
$(X^{(1)},\ldots, X^{(R)})$ is given by
\[
C_{X} = [e^{c_{rr'}\sigma}-1]_{1\le r,r'\le R}.
\]
One easily checks that the mapping $\Psi:[c_{rr'}]\mapsto [e^{c_{rr'}\sigma}-1]$ is not bijective on the set of non negative
symmetric matrices with constant diagonal coefficients: let $ c =\log(1+\kappa)$, $\rho= \log(1+\vartheta)$ where  $\vartheta\ge
-\frac{\kappa}{R-1}$ so that the matrix $[(\kappa-\vartheta)\delta_{rr'} + \vartheta]$ is nonnegative. Assume $\vartheta =
-\frac{\lambda}{R-1}$ with $\lambda \!\in(\log(1+\kappa), \kappa)$. Then this covariance structure cannot be obtained from some
correlated Brownian motions since $c+(R-1)\rho=\log(1+\kappa) +(R-1)\log(1+\vartheta)<0$. 
 
\section{Complexity of the multi-step $R$-$R$ procedure (with consistent Brownian increments) }\label{Complex}
Throughout this section, we assume that the $R$ Gaussian white noises $(U^{(r)}_k)_{1\le k\le nr}$ are consistent
Brownian increments given by~(\ref{ConsistIncr}).
\subsection{Complexity} \label{complexite}The simulation of the correlated sequences $U^{(r)}_k$  can often be neglected in terms
of complexity with respect to the computation of one step of  the $R$ Euler schemes (the main concern being in fact
to spare the random number generator, see below). Then, one can easily derive the complexity of one Monte
Carlo path of this
$R$-order $R$-$R$ extrapolation procedure: for every
$r\!\in\{1,\ldots,R\}$, one has to compute  $r\, n$ values of an Euler scheme. So the complexity $\kappa(R)$ of the
procedure is given by
\[
\kappa(R) = C_{b,\sigma}\times n\times \sum_{r=1}^R r= C_{b,\sigma}\times n\times \frac{R(R+1)}{2}
\]
where $C_{b,\sigma}$ denotes the complexity of a one time step computation in~(\ref{EDS}). Note that using a
recursively iterated  $R$-$R$ procedure would have lead to consider some Euler schemes with steps
$\frac{T}{2^{r-1}n}$,
$1\le r\le R$, (see $e.g.$~\cite{GLA}, p.362 when $R=3$) which in turn would have induced a complexity of order 
$C_{b,\sigma}\times n\times (2^{R}-1)$. 

Thus the global complexity of a Monte Carlo simulation of size $M$ is then
\[
N:= C_{b,\sigma}\times  \frac{R(R+1)}{2}\times n\times M.
\]
Up to a scaling of the global complexity by $1/C_{b,\sigma}$ one may assume without loss of generality that
$C_{b,\sigma}=1$. So evaluating the performances of the  $R$-$R$  extrapolation as a function of its
complexity amounts to solving the following minimization problem:
\[
{\cal E}(R,N):=\min_{ \frac{R(R+1)}{2}\times n\times M=N}\sqrt{\frac{{\rm Var}(f(X_{_T}))}{M}(1+\varepsilon(n))
+
\frac{\widetilde c_{_R}^2}{n^{2R}}(1+O(1/n))}
\]
where $\lim_n \varepsilon(n)=0$, since
\[
\left\|\E(f(X_{_T}))- \frac 1M \sum_{m=1}^M\sum_{r=1}^R
\alpha_r f((\bar X^{(r)}_{T})^m)\right\|_{_2}^2= \frac{{\rm Var}(f(X_{_T}))}{M}(1+\varepsilon(n)) +
\frac{\widetilde c_{_R}^2}{n^{2R}}(1+O(1/n)).
\]
Standard computations lead to the following Proposition.

\begin{Proposition}\label{Complex} Let $R\ge 2$ be the order of expansion of the time discretization error.

\smallskip
\noindent $(a)$ Then,
\[
\lim_{N\to \infty} N^{-\frac12(1-\frac{1}{2R+1})}{\cal E}(R,N)= \frac{1}{\sqrt{2}}\times\Theta(R)
\times R\, |\widetilde c_{_R}|^{\frac{1}{2R+1}}\times \left({\rm Var}(f(X_{_T}))\right)^{\frac{R}{2R+1}}
\]
where $\Theta (R) =  2^{ \frac{1}{2(2R+1)}}R^{-\frac{1}{2R+1}}(1+\frac 1R)^{\frac{R}{2R+1}}\left((2R)^{-
\frac{2R}{2R+1}}+(2R)^{\frac{1}{2R+1}}\right)^{\frac 12}\to 1\;$ as
$\;R\to \infty$.

\medskip
\noindent $(b)$ Furthermore, for a fixed complexity level $N$, the solution $(n(N),M(N))$ of the minimization problem
satisfies
\[
M(N)\sim \frac{2}{R(R+1)}\left(\frac{R+1}{4}\right)^{\frac{1}{2R+1}}\left(\frac{{\rm
Var}(f(X_{_T}))}{\widetilde c_{_R}^2}\right)^{\frac{1}{2R+1}} N^{\frac{2R}{2R+1}}
\]
and
\[ n(N) =\left(\frac{R+1}{4}\right)^{-\frac{1}{2R+1}}\left(\frac{\widetilde c_{_R}^2}{{\rm
Var}(f(X_{_T}))}\right)^{\frac{1}{2R+1}}
  N^{\frac{1}{2R+1}}
\]
as $N\to \infty$ so that
\[
M(N) \sim \frac{1}{2R}\frac{{\rm
Var}(f(X_{_T}))}{\widetilde c_{_R}^2} \,n(N)^{2R}.
\]
\end{Proposition}

\noindent{\sc Comments:}

\smallskip
-- The rate  of convergence increases  with $R$ and tends toward $N^{\frac 12}$. This means that, as
expected,  the higher one expands the time discretization error, the less this error slows down the asymptotic
global rate of convergence.

\smallskip
-- However, an expansion of order $R$ being fixed, the range at which the theoretical rate of convergence
$N^{\frac12(1-\frac{1}{2R+1})}$ does occur depends on the value of  the term
$R  |\widetilde c_{_R}|^{\frac{1}{2R+1}}= \frac{R|c_{_R}|^{\frac{1}{2R+1}}}{(R!)^{\frac{1}{2R+1}}}$ ($\sim \sqrt{eR}|c_{_R}|^{\frac{1}{2R+1}}$
as $R\to \infty$) for which no   bound or estimate is available in practice. Although some theoretical
explicit expression do exist for  the coefficients $c_{_R}$ (at least for small values of $R$, see~\cite{TATU}, \cite{GUY}, \cite{CLKOLA}), it is
usually not possible to have some numerical estimates.

\smallskip
-- The time discretization parameter $n=n(N)$ and the size $M=M(N)$ of the
Monte Carlo simulation are explicit functions of $N$ which involve the
unknown coefficient $|c_{_R}|$. So the above {\em sharp}
$L^2$-rate of convergence is essentially a theoretical bound that cannot
be achieved at a reasonable cost in practice. However one could imagine to produce a rough estimate of $c_{_R}$ by normalizing the
$R$-$R$ extrapolation (by $n^{R}$)  using a small preliminary   MC simulation so as to design  $n(N)$ and $M(N)$ using the above rules.
However the recursive feature of the resulting MC would be lost which seems not very realistic in practice.

\subsection{Efficient simulation of consistent  Brownian increments} The question of
interest   in this section is  the simulation of the Brownian increments. Several methods can be
implemented, the main concern being    to spare the random number generator (its contribution to the
global complexity of the simulation can be neglected in practice as soon as $b$ and $\sigma$ have ``not too
simple" expressions).

The problem amounts to simulating the   increments of the
$R$ Euler schemes, say between absolute time
$t=0$ and $t=T/n$. This means   to simulate  the upper triangular matrix
\[
\left[W_{\frac{\ell T}{rn}}-W_{\frac{(\ell-1) T}{rn}}\right]_{1\le \ell\le r, 1\le r\le R}.
\]
\subsubsection{A lazy approach} Let $R\ge 2$ be a fixed integer. For every $r\!\in\{1,\ldots,R\}$, set
\[
M(R) := {\rm lcm}(1,2,\ldots,R)\qquad \mbox{ and }\qquad m(r) := M(R)/r,\; r=1,\ldots,R.
\]
Then, one simulates $nM(R)$ independent copies $\xi_1,\ldots, \xi_{nM(R)}$  of ${\cal N}(0;1)$ and one sets
\[
U^{(r)}_k = \frac{1}{\sqrt{m(r)}}\sum_{i=(k-1)m(r)+1}^{k m(r)}  \xi_i,\qquad 1\le k\le r\,n.
\]

Such a simulation strategy consumes $n\times M(R)$ random numbers to complete the simulation of the $R$ Euler
schemes until maturity $T$.

\subsubsection{Saving the  random number generator} One only simulates what is needed to get the above matrix. This
means to simulate the Brownian increments between the points of the set
\[
 S_{_R}:= \left\{\frac{\ell}{r}\frac{T}{n},\; 1\le \ell\le r,\; 1\le r\le R\right\}.
\]
Such a simulation strategy consumes $n\times {\rm card}\,S_{_R}$ random numbers to complete the simulation of the $R$
Euler schemes. Obviously
\[
S_{_R}= \left\{\frac{\ell}{r}\frac{T}{n},\; 1\le \ell\le r,\; 1\le r\le R, \; {\rm gcd}(\ell,r)=1\right\}
\]
so that
\[
{\rm card}\,S_{_R}= \sum_{r=1}^R \varphi(r)
\]
where $\varphi$ denotes the Euler function with the convention $\varphi(1)=1$.  Classical Number Theory results  say that
\[
\sum_{r=1}^R \varphi(r)\sim \frac{3}{\pi^2} R^2 \qquad \mbox{ as } \quad R\to +\infty.
\]
whereas $M(R) = e^{(1+o(1))R}$ as $R\to \infty$. In practice, only the first values of
$M(R)$ and card$(S_{_R})$ are of interest. They are
 are reported in the table below.

\bigskip
\centerline{
\begin{tabular}{|c||c|c|c|c|c|c|c|}
\hline
$R$             &1& 2 & 3 & 4  & 5   & 6  & 7\\ \hline\hline
$M(R)$          &1& 2 & 6 & 12 & 60  & 60 & 420 \\ \hline
card$(S_{_R})$&1& 2 & 4 & 6  & 10  & 12 &  18\\ \hline
\end{tabular}
}

\bigskip
For a  given value of $R$, it is necessary to sort the values $\frac{\ell}{r}\!\in
S_{_R}$ and to design the array of resulting lengths in order to simulate the above Brownian matrix at each ``global"
time step $T/n$. But this can be done once since it is universal.

\section{Numerical experiments}\label{Num}

\subsection{Some tables}
We specify below how to simulate the Gaussian white noises $(U_k^{(r)})_{1\le k\le r}$  of the $R$ coupled Euler schemes. 
Note that these tables are equivalent to  a ``multistep" Brownian bridge based recursive simulation of the underlying Wiener
process.

\medskip
-- \fbox{$R=2:$} Let $U_1,\, U_2$ be two i.i.d. copies of ${\cal N}(0;I_q)$. Set
\[
U^{(2)}_i=U_i,\, i=1,2,\qquad U^{(1)}_1= \frac{U_1+U_2}{\sqrt{2}}
\]
\[
\mbox{and }\hskip 5,5 cm \alpha_1=-1,\qquad \alpha_2=2.\hskip 5,5 cm
\]
(Note that $\sum_i\alpha^2_i =5)$).

\medskip
-- \fbox{$R=3:$} Let $U_1,\, U_2,\, U_3,\,U_4$ be four  i.i.d. copies of ${\cal N}(0;I_q)$. Set
\[
U^{(3)}_1=U_1,\quad U^{(3)}_2 = \frac{U_2+U_3}{\sqrt{2}},\quad U^{(3)}_3=U_4,
\]
\[
U^{(2)}_1=\frac{\sqrt{2}\, U_1+ U_2}{\sqrt{3}},\quad U^{(2)}_2 = \frac{U_3+\sqrt{2}\, U_4}{\sqrt{3}},
\]
\[
U^{(1)}_1 =\frac{U^{(2)}_1+U^{(2)}_2}{\sqrt{2}}.
\]
\[
\mbox{and }\hskip 4,5 cm \alpha_1=\frac 12,\qquad \alpha_2= -4,\qquad \alpha_3=\frac 92.\hskip 4,5 cm
\]
(Note that $\sum_i\alpha^2_i =\frac{73}{2}$).

\medskip
-- \fbox{$R=4:$} Let $U_1,\dots, U_6$ be six i.i.d. copies of ${\cal N}(0;I_q)$. Set
\[
U^{(4)}_1=U_1,\quad  U^{(4)}_2 = \frac{U_2+\sqrt{2} \,U_3}{\sqrt{3}},\quad  U^{(4)}_3 = \frac{\sqrt{2} \,U_4+U_5}{\sqrt{3}},\quad
U^{(4)}_4=U_6.
\]
\[
U^{(3)}_1=\frac{ \sqrt{3}\,U_1+ U_2}{2},\quad  U^{(3)}_2 = \frac{U_3+U_4}{\sqrt{2}},\quad
U^{(3)}_3=\frac{U_5+\sqrt{3}\, U_6}{2},
\]
\[
U^{(2)}_1=\frac{U^{(4)}_1+ U^{(4)}_2}{\sqrt{2}},\quad U^{(2)}_2 = \frac{U^{(4)}_3+  U^{(4)}_4}{\sqrt{2}},
\]
\[
U^{(1)}_1 =\frac{U^{(2)}_1+U^{(2)}_2}{\sqrt{2}}.
\]
\[
\mbox{and }\hskip 4 cm \alpha_1=-\frac 16  ,\qquad \alpha_2= 4,\qquad \alpha_3=
-\frac{27}{2},\qquad\alpha_4= \frac{32}{3}.\hskip 4 cm
\]

(Note that $\sum_i\alpha^2_i =\frac{5617}{18}\approx 312.06$).

\medskip

-- \fbox{$R=5:$}
\[
\alpha_1=\frac{1}{24}  ,\qquad \alpha_2= -\frac 83,\qquad \alpha_3=
\frac{81}{4},\qquad
\alpha_4= -\frac{128}{3},\qquad
\alpha_5= \frac{625}{24}.
\]

\subsection{Numerical test with a high volatility Black-Scholes model}\label{BSVanil}

To illustrate the efficiency of the use of consistent white noises in $R$-$R$ 
extrapolation, we  consider the simplest option pricing model, the (risk-neutral) Black-Scholes dynamics, 
but with unusually high volatility. (We are aware this model used to price Call options does not fulfill the theoretical assumptions made above).
To be precise
\[
dX_t= X_t\,(rdt+\sigma dW_t), 
\]
with the following values for the
parameters

%\bigskip
%\centerline{
%\fbox{%$\displaystyle
\[
X_0=100, \; K=100, \; r= 0.15,\; \sigma = 1.0, \; T=1.
\]
%$.}}
%
%\bigskip
Note that a volatility $\sigma= 100\%$ per year is equivalent to a 4 year maturity with volatility $50\%$ (or 16 years with
volatility $25 \%$).  The reference Black-Scholes premium is $C^{BS}_0= 42.96$. We consider the Euler scheme with step
$T/n$ of this equation.
\[
\bar X_{t_{k+1}}= \bar X_{t_k}\left(1+r\frac Tn +\sigma \sqrt{\frac Tn}U_{k+1}\right),\qquad \bar X_{0}=X_0,
\]
where $t_k= \frac{kT}{n}$, $k=0,\ldots,n$. We want to price a vanilla Call option  $i.e.$ to compute
\[
C_0 = e^{-rT} \E((X_{_T}-K)_+)
\]
using a Monte Carlo simulation with $M$ sample paths.

\medskip 
\noindent $\bullet$ {\sc Three step $R$-$R$ extrapolation ($R=3$):} We processed three Monte Carlo simulations of common size $M=10^6$
to evaluate the efficiency of the $R$-$R$ extrapolation with $R=3$, having in mind that, for a given complexity $N$, teh
size of teh Monte Carlo simulation and the time discritization parameter satisfy $M(N)\propto n(N)^{2R}$, see
Proposition~\ref{Complex}$(b)$).

\smallskip
-- An $R$-$R$ extrapolation of order $R=3$ as defined
by~(\ref{RmultiRomb})  and~(\ref{alpha}) with $n= 2,\,4,\,6,\,8,\, 10$
with {\em consistent increments} (the maximal number of steps is $R\,n=3\,n$). Note these specifications for $n$ are quite low
in comparison with the high volatility of the $B$-$S$ model. 

\smallskip
-- An $R$-$R$ extrapolation of order $R=3$ with the same architecture but implemented with {\em independent increments}.

\smallskip
-- A {\em regular  Euler scheme with the same complexity} $i.e.$ with
$R(R+1)/2\times n\!=\! 6\,n$ steps.

\medskip The results are depicted in the Figures~1 and~2 below. In Figure~1, the abscissas represent the size ($3n$) of the
Euler scheme with the highest discretization frequency used in the procedure. In Figure~2, the abscissas represent the
size ($6n$)  of the standard Euler scheme with the
same complexity as the   $R$-$R$ extrapolation ($R=3$). The main conclusions are   the following:

-- The standard deviation of the $R$-$R$ extrapolation with independent noises is $5$ times greater than the one observed with consistent increments. This  
makes  this higher order $R$-$R$ expansions ($R=3$) useless at least within the usual range of our Monte Carlo
simulations (see~Fig.~1): in fact it is less efficient in our high volatility setting than the Euler scheme with equivalent complexity and
less efficient than a standard $R$-$R$ extrapolation ($R=2$).

--  Considering consistent increments in the $R$-$R$ extrapolation gives
 the method its full efficiency as emphasized by Figure~2 : $R$-$R$ extrapolation is clearly much  more efficient than the Euler
scheme of equal complexity in a Monte Carlo simulation of size $M=10^6$.
\begin{figure}
\begin{center}
\begin{tabular}{cc}
$$\includegraphics[height =5cm,width=7.5cm]{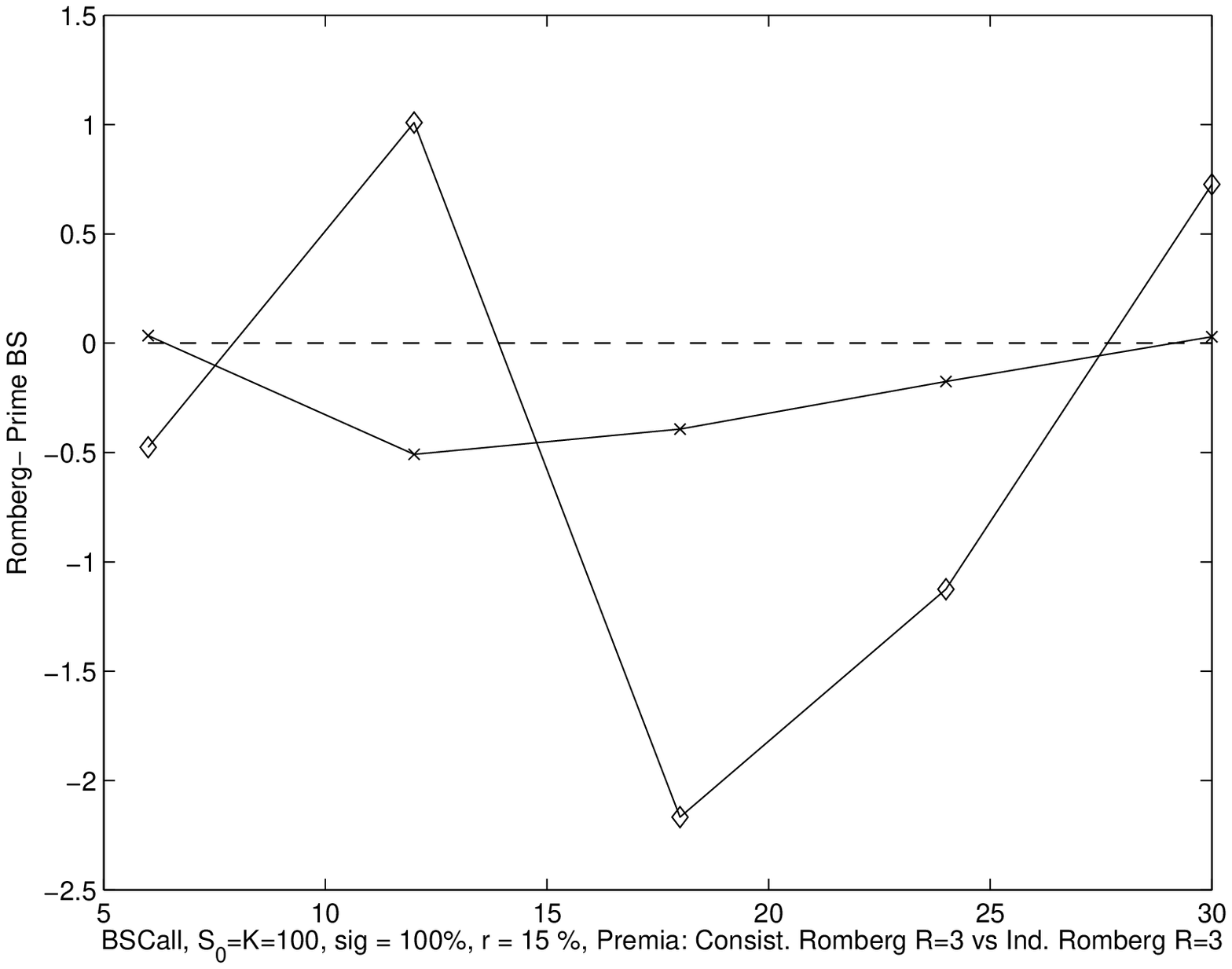}$$&
$$\includegraphics[height =5cm,width=7.5cm]{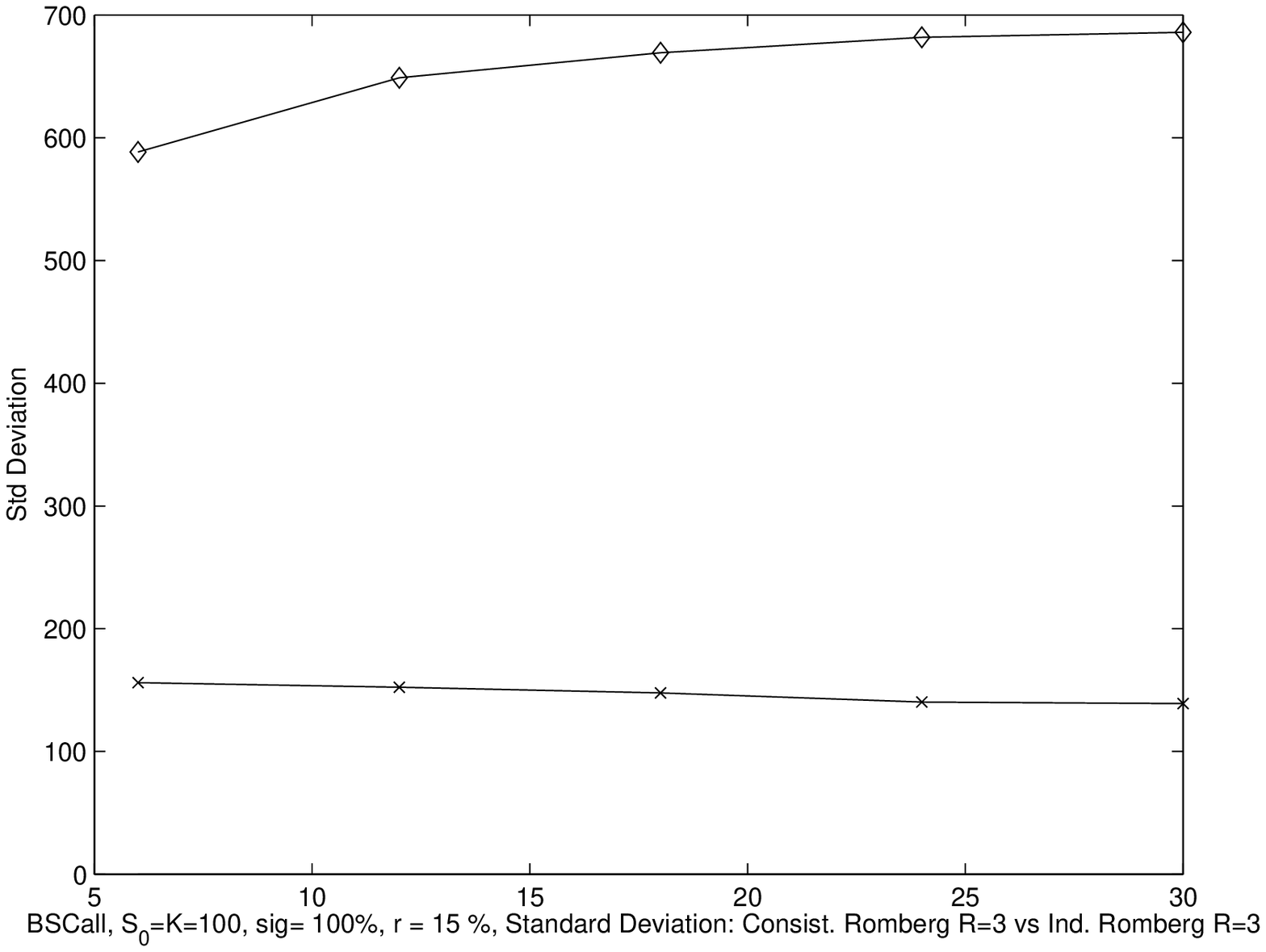}$$
%&$$\includegraphics[height =4 cm,width=4cm]{Fig3(erg).eps}$$
 \end{tabular}
\caption{\small{\sc   B-S Euro Call option.} $M\!=\!10^6$. {\em $R$-$R$ extrapolations $R\!=\!3$.  Consistent Brownian increments
($\!\times\!\!-\!\!\!-\!\!\times\!\!-\!\!\!-\!\!\times\!$)
 $vs$ independent Brownian increments ($\diamondsuit\!$---$\diamondsuit\!$---$\diamondsuit\!$).
 $X_0\!=\!K\!=\!100$,
$\sigma\!=\!100 \%$, $r\!=\!15 \%$. Abscissas: $3n$, $n=2,4,6,8,10$.
 {\rm Left:} Premia. {\rm Right:} Standard Deviations.}}
\label{fig:CallBS_RR}
\end{center}
\end{figure}
\begin{figure}
\begin{center}
\begin{tabular}{cc}
$$\includegraphics[height =5cm,width=7.5cm]{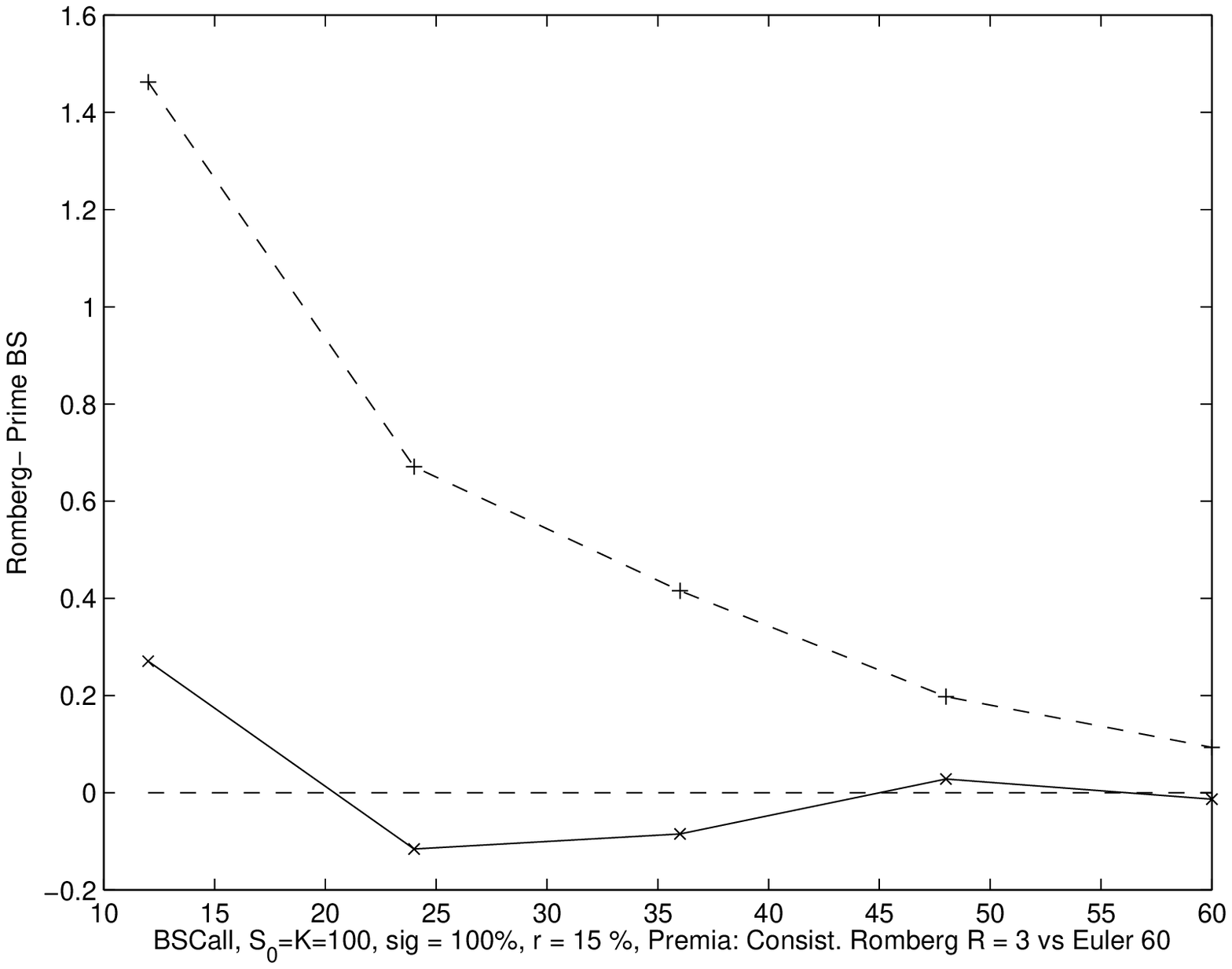}$$&
$$\includegraphics[height =5cm,width=7.5cm]{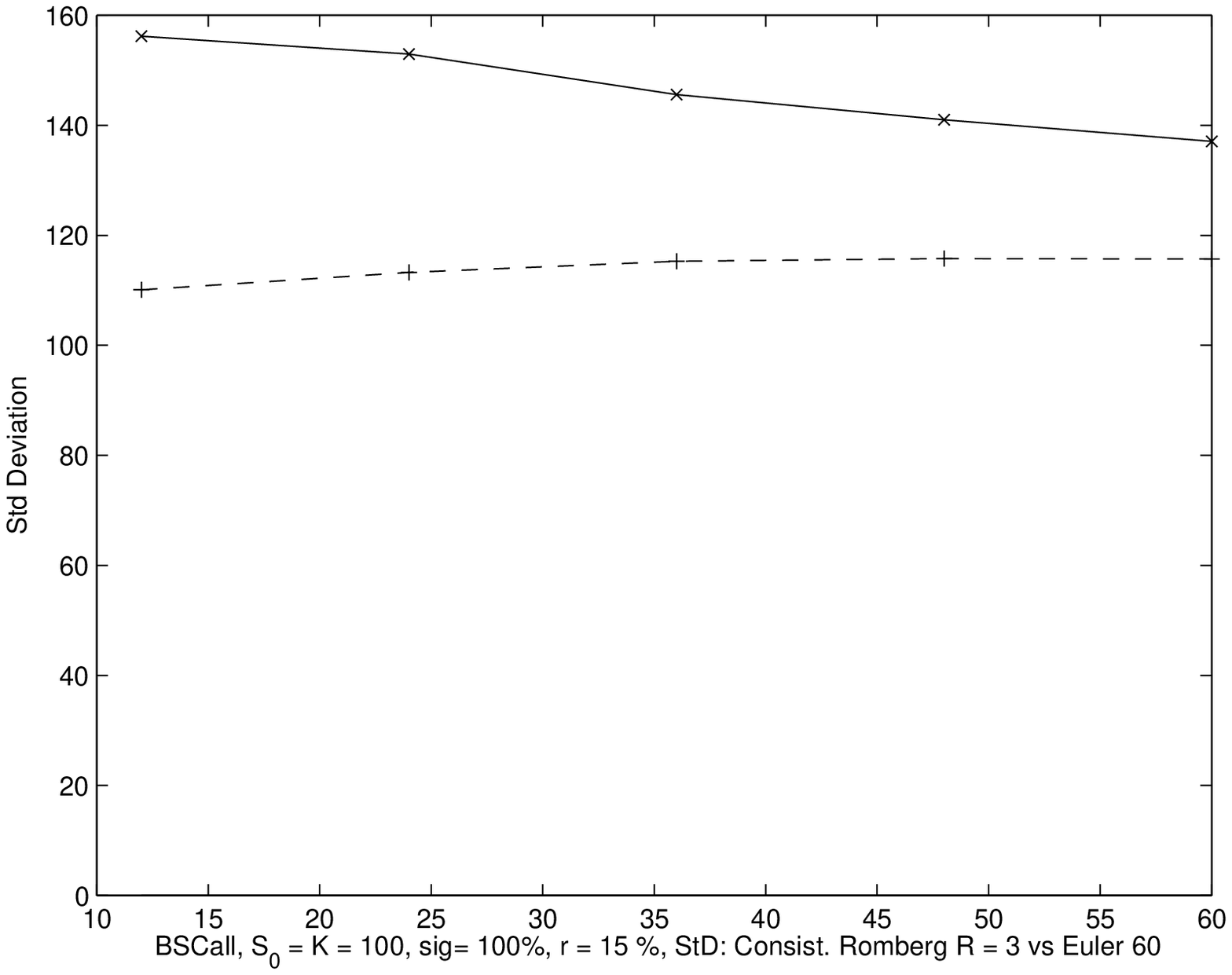}$$
%&$$\includegraphics[height =4 cm,width=4cm]{Fig3(erg).eps}$$
 \end{tabular}
\caption{\small{\sc   B-S Euro Call option.} $M\!=\!10^6$. {\em  $R$-$R$ 
extrapolation $R\!=\!3$. Consistent Brownian
increments ($\!\times\!\!-\!\!\!-\!\!\times\!\!-\!\!\!-\!\!\times\!$) $vs$ Euler scheme with equivalent complexity
($\!+\!$ - - $\!+\!$ - - $\!+\!$). $X_0\!=\!K\!=\!100$, $\sigma\!=\!100 \%$, $r\!=\!15
\%$. Abscissas: $3n$, $n=2,4,6,8,10$. Abscissas: $6n$, $n=2,4,6,8,10$.
 {\rm Left:} Premia. {\rm Right:} Standard
Deviations.}}
\end{center}
\end{figure}

\medskip 
\noindent $\bullet$ {\sc Four step $R$-$R$ extrapolation ($R\!=\!4$):} With  a size   of $M=10^6$, the $R$-$R$-extrapolation with $R=4$
is not completely convincing: at this range, the variance of the estimator is not  yet controlled by ${\rm Var}(f(X_{_T}))$ for  
the selected (small)  values of the discretization parameter $n$. 

However, for  a larger simulation, say $M=10^8$, the multistep extrapolation of order $R=4$ clearly  becomes the most efficient one as
illustrated by Figure~3 (right) and Table~1.

To carry out a comparision, we  implemented this time:

\smallskip
-- Two $R$-$R$ extrapolations with  orders $R=3$ and $R=4$ as defined by~(\ref{RmultiRomb}) and~(\ref{alpha}) with $n =
2,\,4,\,6,\,8,\, 10$ with {\em consistent increments}. 

\smallskip
-- A {\em regular  Euler scheme with the same complexity} $i.e.$ with
$4\times 5/2\times n\!=\! 10\,n$ steps.

\medskip
In Figure~3, the abscissas represent the size ($10n$) of the Euler scheme with same complexity as the $R$-$R$ extrapolation with $R=4$.

\bigskip
\begin{tabular}{c||c c c c c }
%\hline
$n$   & 2&4&6&8&10\\
\hline
\hline
{\bf BS} $Ref.$ &  & &$\mathbf{42.96}$& &  \\
\hline
$R=3$& $42.93$ &  $42.55$ &  $42.80$  &  $42.90$ &  $42.95$\\
					& ($-0.07\%$)&   ($-0.94\%$)& ($-0.37\%)$ &  ($-0.14\%$) &  ($-0.01\%$)\\
$R=4$&Ê$42.28$ &  $42.92$ &  $42.97$  & $42.94$  & $42.97$\\
      &($-1.59\%$) & ($-0.08\%$) & ($0.04\%$) & ($-0.03\%$)  & ($0.03\%$)
%\hline
\end{tabular}

%Ref BS  42.9571

%$R=3$  42.9268   42.5537   42.7976   42.8965   42.9521
%  155.5661  151.4581  145.3872  140.7010  137.3198
%
%$R=4$ 42.2764   42.9246   42.9737   42.9455   42.9692
%  326.6174  293.7677  262.8094  240.6565  224.8540
  
%\bigskip
\begin{figure}\label{fig3}
\begin{center}
\begin{tabular}{cc}
$$
\includegraphics[height =5cm,width=7.5cm]{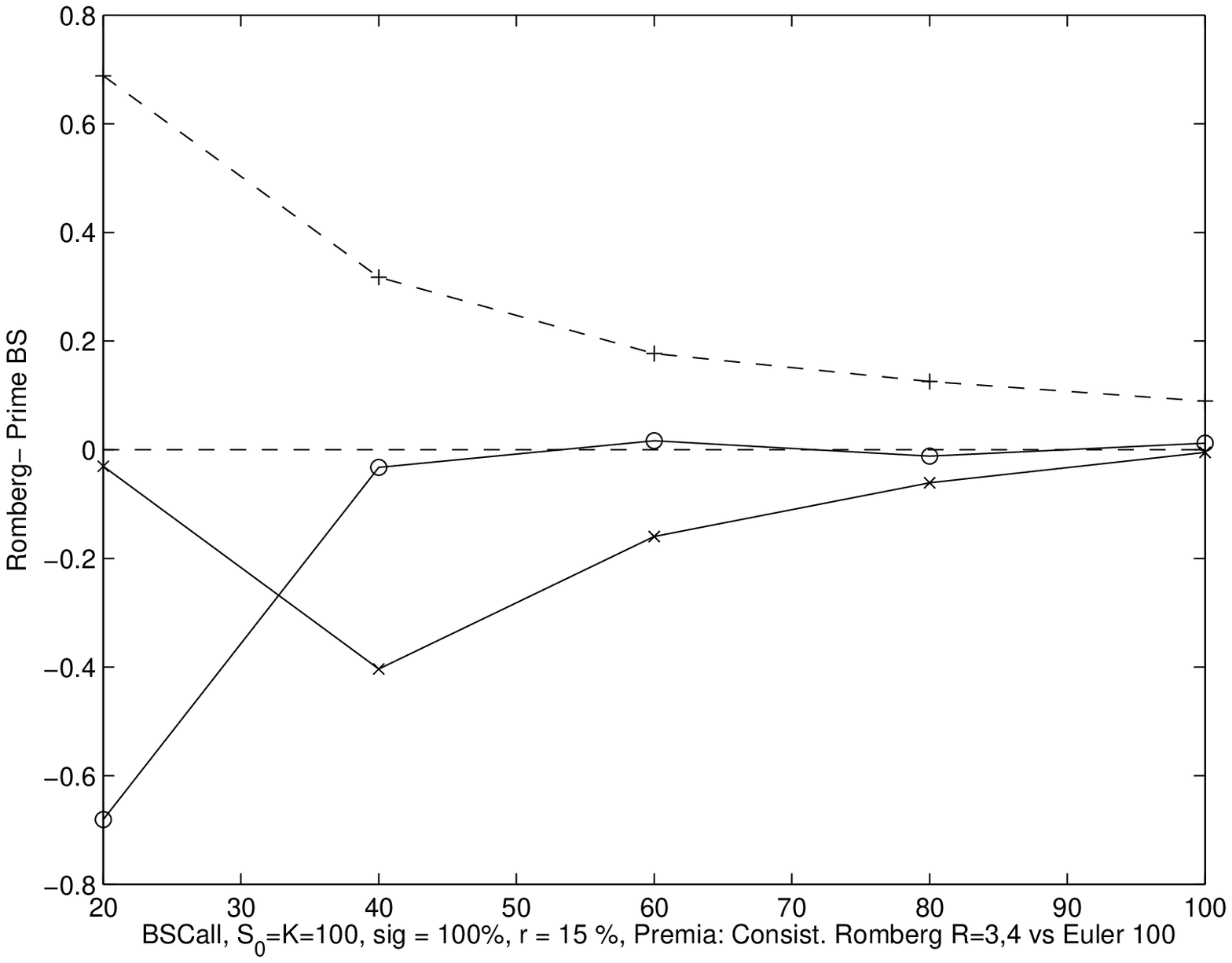}$$&
$$
\includegraphics[height =5cm,width=7.5cm]{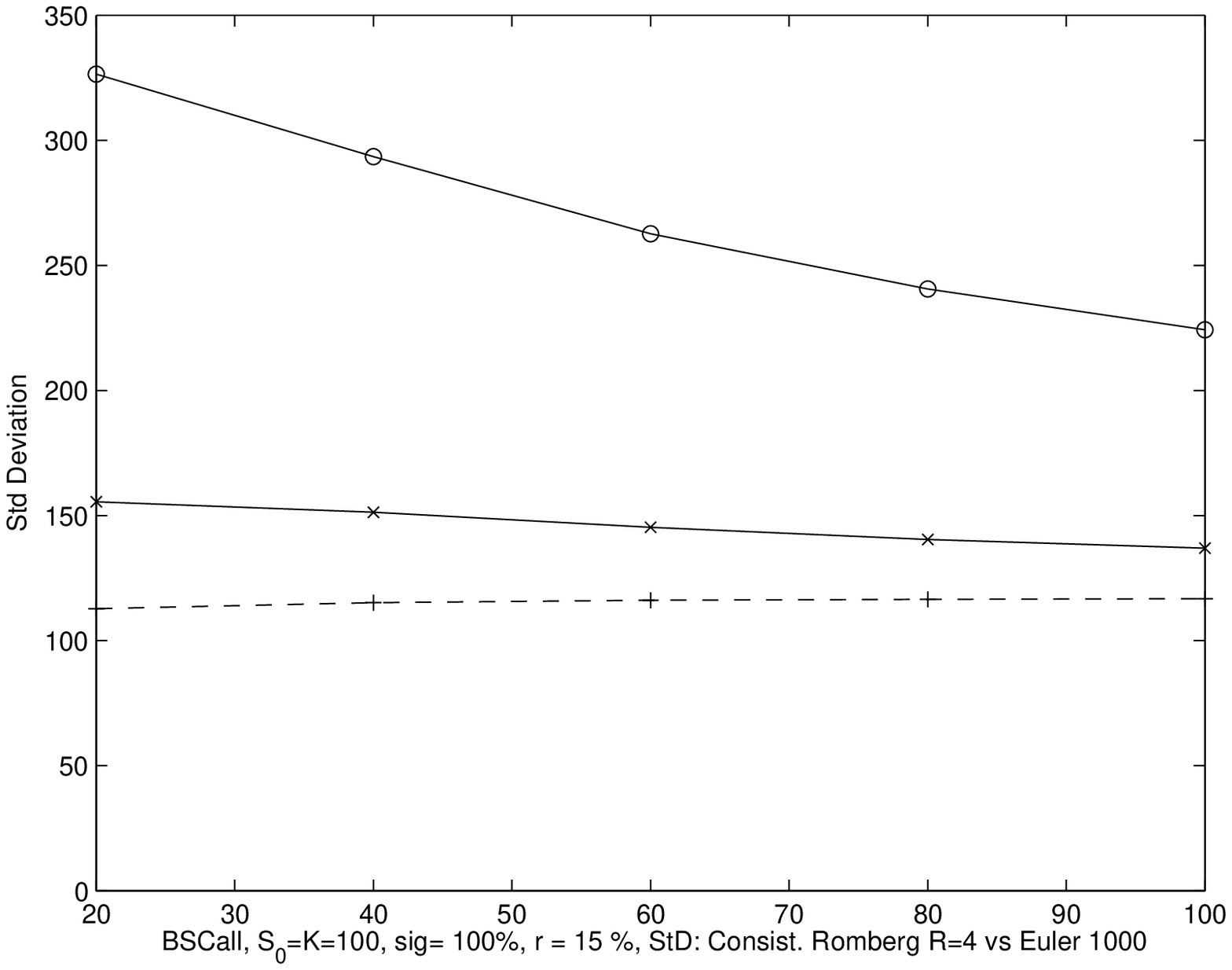}
$$
 \end{tabular}
\caption{\small{\sc   B-S Euro Call option.}  $M\!=\!10^8$. {\em  $R$-$R$ 
extrapolation $R\!=\!3,\, 4$. Consistent Brownian
increments $R=3$: ---$\!\times\!$---$\!\times\!$---$\!\times\!$---; Consistent Brownian
increments $R=4$: $-\!\!\!-\!o\!-\!\!\!-\!o\!-\!\!\!-\!o\!-\!\!\!-$;  Euler scheme with equivalent complexity
($\!+\!$ - - $\!+\!$ - - $\!+\!$). $X_0\!=\!K\!=\!100$, $\sigma\!=\!100 \%$, $r\!=\!15
\%$. Abscissas: $10n$, $n=2,4,6,8,10$.
 {\rm Left:} Premia. {\rm Right:} Standard
Deviations.}}
\end{center}
\end{figure}

\medskip
\centerline{Table~1. {\em $R$-$R$ extrapolation with $R=3, \,4$ $vs$ the Euler scheme with step
$1/(10n)$.}} 
\medskip 
This emphasizes   the natural field of application of multistep $R$-$R$ extrapolation for numerical applications ($R\ge 4$): this
is the most efficient method to obtain  accurate results in a high variance framework: it allows a smaller size of discretization
step.

\subsection{Further numerical experiments:  path dependent options}\label{BSpathdep}

In this section we will consider some path-dependent (European) options $i.e.$ related to some payoffs $F((X_t)_{t\in[0,T]})$ where $F$
is a functional defined on the set  
$\D([0,T],\R^d)$   of right continuous left-limited functions $x:[0,T] \to \R$. It is clear that
all the  asymptotic control of the variance  obtained in Section~\ref{SecCont} for the estimator $\sum_{r=1}^R \alpha_r f(\bar X^{(r)}_{_T})$ of $\E(f(X_{_T}))$ when $f$ is continuous can be extended to functionals
$F:\D([0,T],\R^d)\to \R$  which are  $\P_{X}$-$a.s.$ continuous  with respect to the $\sup$-norm defined by
$\|\,x\,\|_{_{\sup}}:=\sup_{t\in [0,T]}|x(t)|$ with polynomial growth ($i.e.$
$|F(x)| =O(\|x\|^\ell_{_{\sup}})$ for somme natural integer $\ell$ as
$\|x\|_{_{\sup}}\to\infty$). This simply follows
from the fact that the (piecewise constant) Euler scheme
$\bar X$ (with step $T/n$)  converges for the
$\sup$-norm toward $X$ in $L^2(\P)$. 

The same result holds with the continuous Euler scheme $\bar X^c$  defined by 
$$
\forall\, t\!\in [0, T],\qquad \bar X^c_t = x_0+
\int_0^tb(\bar X^c_{\underline s})ds +\int_0^t\sigma(\bar X^c_{\underline s})dW_s, \;\underline t=\lfloor nt/T\rfloor
$$ 
(the diffusion coefficients have been frozen between time discretization instants). With this scheme, one can  simply consider the path set ${\cal
C}([0,T],\R^d)$ instead of  $\D([0,T],\R^d)$.

Furthermore, this  asymptotic control of the variance  holds true  with any $R$-tuple $\alpha=(\alpha_r)_{1\le r\le
R}$ of weights coefficents satisfying $\sum_{1\le r\le R}\alpha_r =1$, so these coefficients can be adapted to
the structure of the weak error expansion. 

On the other hand, in the recent past years, several papers provided  some weak rates of convergence for some families of
functionals $F$. These works were essentially motivated by the pricing of path-dependent (European) options, like Asian, lookback
or barrier options. This corresponds to functionals 
\[
F(x) := \Phi(\int_0^T x(s)ds),\quad F(x) := \Phi(x(T), \sup_{t\in [0,T]}x(t), \inf_{t\in [0,T]}x(t)), \quad F(x)=
\Phi(x(T))\mbox{\bf 1}_{\{\tau_{_D}(x)\le T\}}
\]
where $\Phi$  is usually  at least Lipschitz  and $\tau_{_D}:=\inf\{s\!\in[0,T],\, x(s\pm)\!\in\, ^c\!D\}$ is the hitting time of $^c\!D$ by
$x$ (\footnote{when $x$ is stepwise constant and c\`adl\`ag, one can  write ``$s$" instead of  "$s\pm$".}). Let us briefly mention two well-known examples:

\smallskip

-- In~\cite{GOB}, it is established that if the domain $D$ has a smooth enough boundary, $b,\,\sigma\!\in {\cal C}^3(\R^d)$,
$\sigma$ uniformly elliptic on $D$, then for every Borel bounded function $f$ vanishing in a neighbourhood of
$\partial D$,
\begin{equation}\label{Bar0}
\E(f(\bar X_{_T})\mbox{\bf 1}_{\{\tau(\bar X)>T\}}) - \E(f(X_{_T})\mbox{\bf
1}_{\{\tau(X)>T\}})=O\left(\frac{1}{\sqrt{n}}\right)\quad \mbox{ as }\quad n\to \infty.
\end{equation}
If furthermore, $b$ and $\sigma$ are ${\cal C}^5$, then
\begin{equation}\label{Barc}
\E(f(\bar X^c_{_T})\mbox{\bf 1}_{\{\tau(\bar X^c)>T\}}) - \E(f(X_{_T})\mbox{\bf
1}_{\{\tau(X)>T\}})=O\left(\frac{1}{n}\right)\quad \mbox{ as }\quad n\to \infty.
\end{equation}
Note however that these assumptions are not satisfied by usual barrier options (see below). 

\smallskip
--  it is suggested  in~\cite{SETO} (including a rigorous proof when $X=W$)  that if $b,\, \sigma\!\in {\cal C}^4_b(\R)$,
$\sigma$ is uniformly elliptic and
$\Phi\!\in{\cal C}^{4,2}(\R^2)$ (with some  partial derivatives with polynomial  growth), then
\begin{equation}\label{Lkb0}
\E(\Phi(\bar X_{_T}, \min_{0\le k\le n}\bar X_{t_k}) - \E(\Phi(X_{_T},\min_{t\in [0,T]}X_t)) =
O\left(\frac{1}{\sqrt{n}}\right)\quad \mbox{ as }\quad n\to \infty.
\end{equation}
A similar improvement -- $O(\frac 1n)$ rate --  as above can be expected when replacing $\bar X$  by   the continuous
Euler scheme $\bar X^c$.

\medskip For both classes of functionals (with  $D$ as a half-line in $1$-dimension in the first   setting), the practical implementation of the continuous Euler scheme is known as the {\em Brownian bridge method}. It relies on the simulations of  the distribution of
$\min_{t\in [0,T]} \bar X^c _t$ and $\max_{t\in [0,T]}
\bar X^c _t$ given the event $\{\bar X^c_{t_k}=x_k,\; k=0,\ldots,n\}$ ($=\{\bar X_{t_k}=x_k,\; k=0,\ldots,n\}$).  This distribution is known since
\[
{\cal L}(\max_{t\in [0,T]}
\bar X^c _t\,|\, \{\bar X_{t_k}=x_k,\; k=0,\ldots,n\}) = {\cal L}(\max_{0\le k\le n-1}G^{-1}_{x_k,x_{k+1}}(U_k))
\]
where 
$$
G^{-1}_{x,y}(u)= \frac 12\left(x+y+\sqrt{(y-x)^2-2T\sigma^2(x)\log(u)/n}\right)
$$ 
and $(U_k)_{0\le k\le n-1}$ are $i.i.d.$
uniformly distributed random variables over the unit interval. A similar formula holds for the minimum using now  the
inverse distribution function 
$$
F^{-1}_{x,y}(u)= \frac 12\left(x+y-\sqrt{(y-x)^2-2T\sigma^2(x)\log(u)/n}\right).
$$

At this stage there are two   ways  to implement the (multistep)
$R$-$R$ extrapolation with consistent Brownian increments in order to improve the performances of the original  (stepwise constant or
continuous) Euler schemes. Both rely on natural conjectures about the existence of a higher order expansion of the time discretization
error suggested by the above rates of convergence~(\ref{Bar0}),~(\ref{Barc}) and~(\ref{Lkb0}). 

\medskip
$\bullet$ {\em Standard Euler scheme:} As concerns the standard Euler scheme, this means the
existence of a vector space $V$ (stable by product) of admissible functionals satisfying
\begin{equation}\label{Expansion2}
({\cal E}^{{\frac 12},V}_{_R})\hskip 0,75 cm \equiv \hskip 0,75 cm\forall\, F\!\in V,\qquad  \E(F(X))= \E(F(\bar
X)) + \sum_{k=1}^{R-1} \frac{c_k}{n^{\frac k2}} + O(n^{-\frac{R}{2}}).\hskip 0,75 cm 
\end{equation}
 
 The main point for practical application is to compute the weights $\alpha^{(\frac 12)}= (\alpha^{(\frac 12)}_r)_{1\le
r\le R}$ of the extrapolation are modified. Namely
\[
\alpha^{(\frac 12)}_r= \alpha^{(\frac 12)}(r,R) := \frac{(-1)^{R-r}}{2} \frac{r^R}{r!(R-r)!} \prod_{k=1}^R\left(1+\sqrt{\frac k
r}\right),\qquad 1\le r\le R.
\]

For small values of $R$, we have
\begin{eqnarray*}
\mbox{\fbox{$R=2$}}\!&\!\!&\! \alpha^{(\frac 12)}_1 = -(1+\sqrt{2}),\qquad \alpha^{(\frac 12)}_2= \sqrt{2}(1+\sqrt{2}).\\
*[.4em] \mbox{\fbox{$R=3$}}\!&\!\!&\!\alpha^{(\frac 12)}_1 = \frac{\sqrt{3}-\sqrt{2}}{2\sqrt{2}-\sqrt{3}-1},\qquad
\alpha^{(\frac 12)}_2=-2\frac{\sqrt{3}-1}{2\sqrt{2}-\sqrt{3}-1},\qquad \alpha^{(\frac 12)}_3
= 3\frac{\sqrt{2}-1}{2\sqrt{2}-\sqrt{3}-1}.\\
*[.4em] \mbox{\fbox{$R=4$}}\!&\!\!&\!\alpha^{(\frac 12)}_1 = -\frac{(1\!+\!\sqrt{2})(1\!+\!\sqrt{3})}{2},\qquad
\alpha^{(\frac 12)}_2=4(\frac 3 2\!+\!\sqrt{2})(\sqrt{3}\!+\!\sqrt{2}),\\ 
*[.4em]\!&\!\!&\!\alpha^{(\frac 12)}_3= -\frac 3 2(\sqrt{3}\!+\!\sqrt{2})(2\!+\!\sqrt{3})(3\!+\!\sqrt{3}),\qquad\alpha^{(\frac
12)}_4= 4(2\!+\!\sqrt{2})(2\!+\!\sqrt{3}).
\end{eqnarray*}

Note that these coefficients have greater absolute values than in the standard case. Thus if $R=4$, $\sum_{1\le r\le
4}(\alpha^{(\frac 12)}_r)^2\approx 10\, 900$! which induces an increase of the variance term for too small values of the time
discretization parameters $n$ even when  increments are consistently generated.  The complexity computations carried out in
Section~\ref{complexite} need to be updated but {\em grosso modo} the optimal choice for the  time discretization parameter $n$
as a function of the $MC$ size
$M$ is 
\[
   M\propto n^{R}. 
\]

$\bullet$ {\em The continuous Euler scheme:} The conjecture is simply to assume that the expansion  $({\cal E}^{V}_{_R})$ now
holds for a vector space $V$ of functionals $F$ (with linear growth with respect to the sup-norm). The  increase of the
complexity induced by the Brownian bridge method  is difficult to quantize: it amounts  to computing $\log(U_k)$ and the inverse
distribution functions $F^{-1}_{x,y}$ and $G^{-1}_{x,y}$.

The second difficulty is that  simulating  (the extrema of) some of  continuous Euler schemes using the Brownian
bridge  in a consistent way is not straightforward at all. However, one can reasonably expect that using independent Brownian
bridges ``relying" on stepwise constant Euler schemes with consistent Brownian increments will have a small impact on the global
variance (although slightly increasing it). 

\bigskip

To illustrate and compare these   approaches we carried some numerical tests on partial
lookback and barrier options  in the Black-Scholes model presented in the previous section. 

\bigskip
\ni $\rhd$ {\sc Partial lookback options:} The partial lookback Call option is defined by its payoff functional
\[
F(x)= e^{-rT}\left(x(T)-\lambda\min_{s\in [0,T]}x(s)\right)_+, \quad x\!\in {\cal C}([0,T],\R),
\]
where $\lambda>0$ (if $\lambda\le 1$, the $(\,.\,)_+$ can be dropped).  The premium 
\[
{\rm Call}^{Lkb}_0 = e^{-rT} \E((X_{_T}-\lambda\min_{t\in [0,T]}X_t)_+)
\]
 is
given by
\[
 {\rm Call}^{Lkb}_0 = X_0{\rm Call}^{BS}\left(1,\lambda , \sigma, r, T\right)+\lambda\frac{\sigma^2}{2r} X_0 {\rm
Put}^{BS}\left(\lambda^{\frac{2r}{\sigma^2}}, 1, \frac{2r}{\sigma},r,T\right).
\]
We took the same values for the $B$-$S$ parameters as in the  former section and set the coefficient $\lambda$ at $\lambda =1.1$.
For this set of parameters
${\rm Call}^{Lkb}_0= 57.475$.
%57.4746

 As
concerns the MC simulation size, we still set $M=10^6$. We compared the following three methods for every choice of $n$:
 
\smallskip
--  A $3$-step $R$-$R$ extrapolation ($R=3$) of the stepwise constant Euler scheme (for which a $O(n^{-\frac 32})$-rate can be 
expected from the conjecture).

\smallskip
-- A $3$-step $R$-$R$ extrapolation ($R=3$) based on the continuous Euler scheme(Brownian bridge method) for which a
$O(\frac{1}{n^3})$-rate can be conjectured (see~\cite{GOB}).

\smallskip
--   A continuous  Euler scheme  (Brownian bridge method) of equivalent complexity $i.e.$ with discretization parameter $6n$ for which
a $O(\frac 1n)$-rate can be expected (see~\cite{GOB}).

\smallskip
The  three procedures have the same complexity if one neglects the cost of the bridge simulation with respect to that
of the diffusion coefficients (note this is  very conservative in favour of ``bridged schemes").

\medskip We do not reproduce the results obtained for the standard stepwise constant Euler scheme which are clearly out of the game
(as already emphasized in~\cite{GOB}).  In Figure~4, the abscissas represent the size of Euler scheme with equivalent complexity
($i.e.$ $6n$, $n=2,4,6,8,10$). Figure~4$(a)$ (left) shows that both $3$-step $R$-$R$ extrapolation methods  converge
significantly faster than the   ``bridged" Euler scheme with equivalent complexity in this high volatility framework. The standard
deviations depicted in Figure~4(a) (right) show that the $3$-step $R$-$R$ extrapolation of the Brownian bridge is controlled even for
small values of $n$. This  is not the case with the $3$-step $R$-$R$ extrapolation method of the stepwise constant Euler scheme. Other
simulations --~not reproduced here~-- show this is already true for the standard
$R$-$R$ extrapolation and the bridged Euler scheme. In any case the multistep
$R$-$R$ extrapolation with $R=3$ significantly outperforms the bridged Euler scheme.

When $M=10^8$, one verifies (see Figure~$4(b)$) that the time discretization error of the $3$-step $R$-$R$ extrapolation   vanishes like
for the partial lookback option. In fact for $n=10$ the $3$-step bridged Euler scheme yields a premium equal to $57.480$ which
corresponds to less than   half a cent error, $i.e.$  $0.05\,\%$ accuracy! This result being obtained without any control variate
variable. 

The $R$-$R$ extrapolation of the standard Euler scheme also provides excellent results. In fact it seems difficult to discriminate
them with those obtained with the bridged schemes, which is slightly unexpected if one think about the natural conjecture about the
time discretization error expansion.

\medskip
As a theoretical conclusion, these results strongly support both conjectures about the  existence of expansion for  the weak error in
the $(n^{-p/2})_{p\ge 1}$ and $(n^{-p})_{p\ge 1}$ scales respectively.

\begin{figure}
\begin{center}
\begin{tabular}{lcc}
(a)&&\\ $\includegraphics[height =5cm,width=7.5cm]{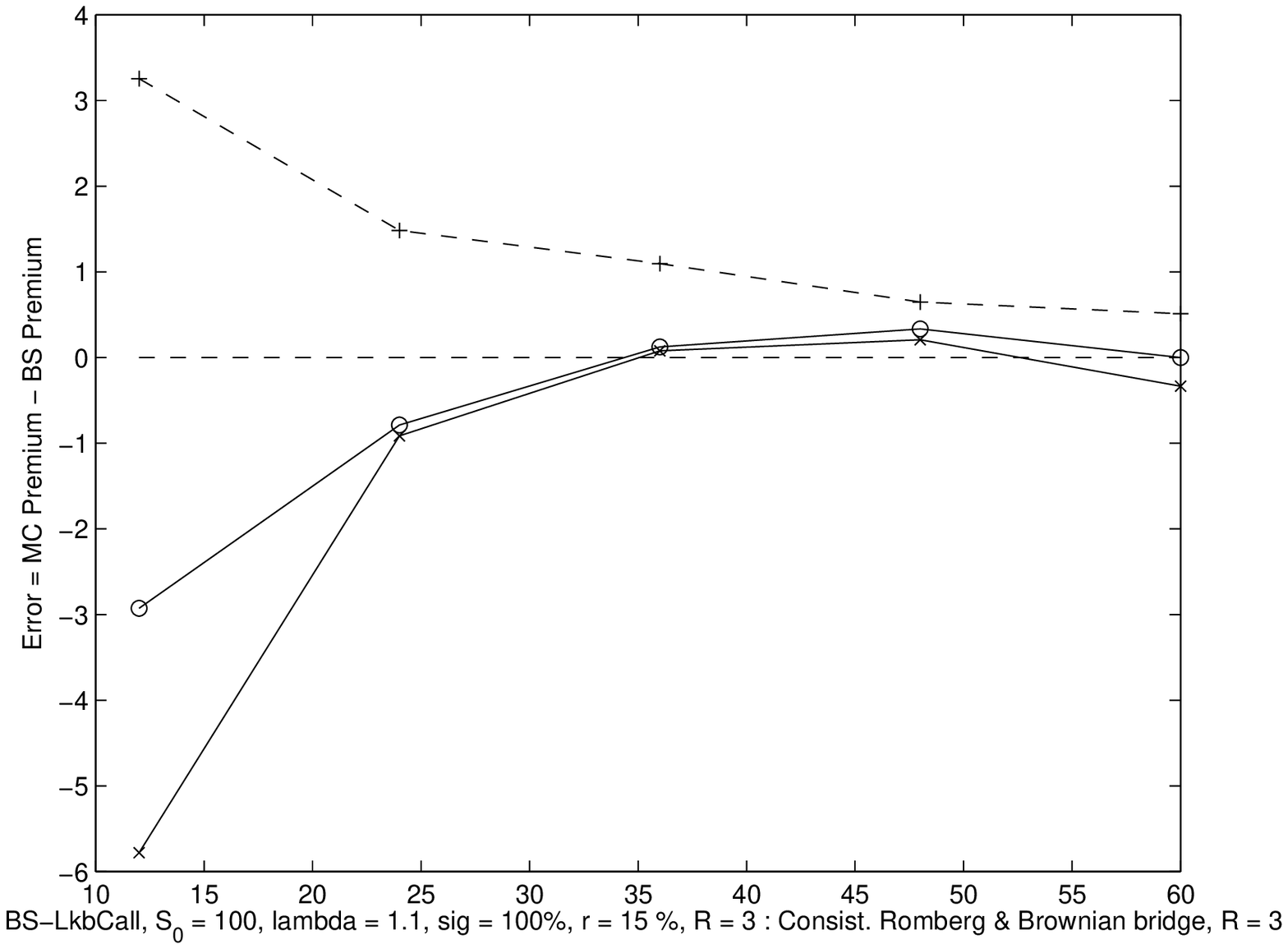}$&
$\includegraphics[height =5cm,width=7.5cm]{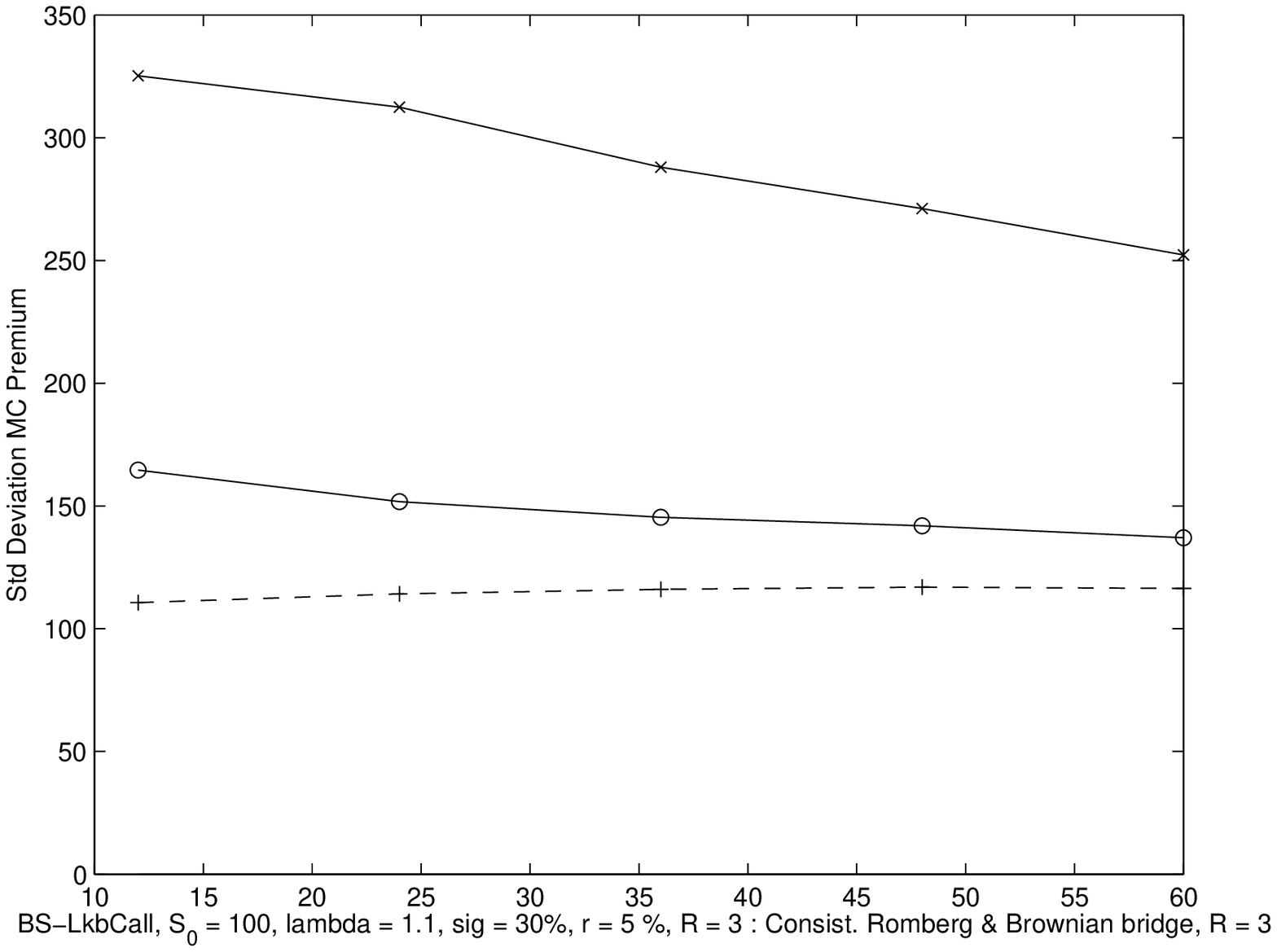}$ \\
(b)&&\\ $\includegraphics[height =5cm,width=7.5cm]{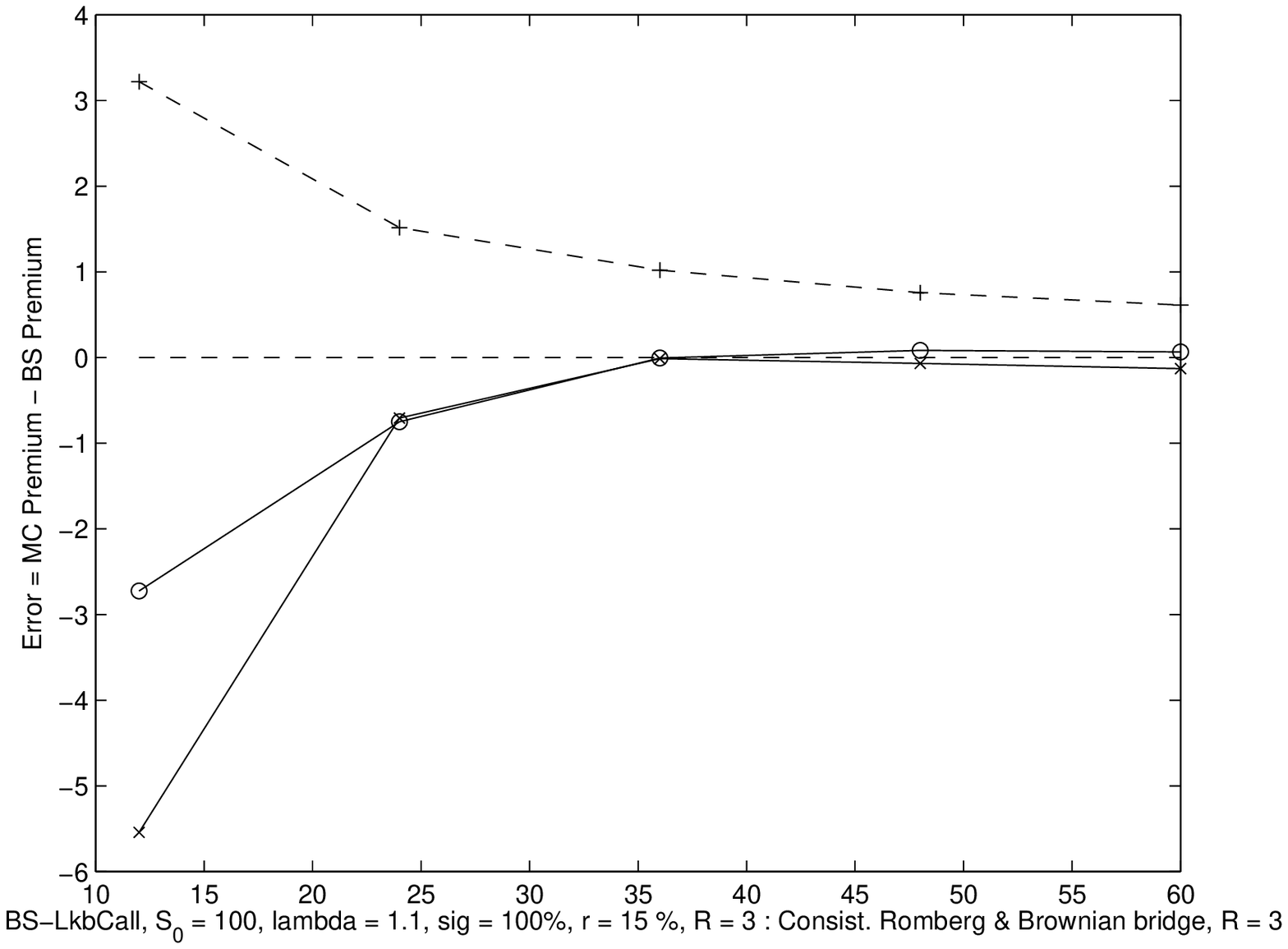}$&
 $\includegraphics[height =5cm,width=7.5cm]{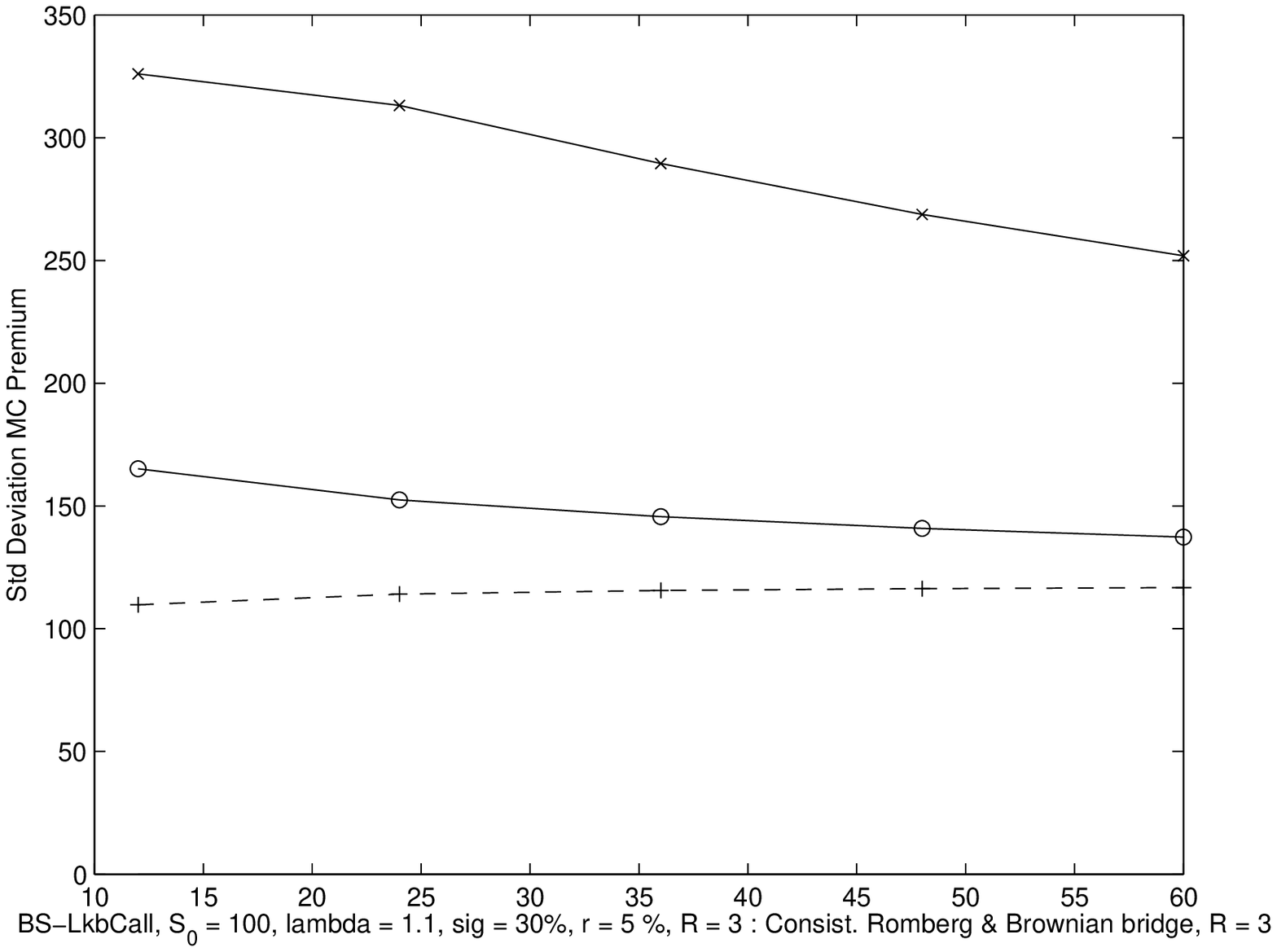}$
%&$$\includegraphics[height =4 cm,width=4cm]{Fig3(erg).eps}$$
 \end{tabular}
\caption{\small {\sc $B$-$S$  Euro  Partial Lookback Call option.} $(a)$ $M\!=\!10^6$.  {\em  $R$-$R$ extrapolation  ($R\!=\!3$) of the
Euler scheme with Brownian bridge: $o\!-\!\!\!-\!o\!-\!\!\!-\!o$.   Consistent  $R$-$R$ 
extrapolation  ($R\!=\!3$): $\times\!\!$---$\!\!\times\!\!$---$\!\!\times$. Euler scheme with Brownian bridge with equivalent
complexity: $+--+--+$. $X_0\!=\! 100$,  $\sigma\!=\!100 \%$, $r\!=\!15\,\%$, $\lambda =1.1$. Abscissas: $6n$, $n=2,4,6,8,10$.
 {\rm Left:} Premia. {\rm Right:} Standard Deviations.} $(b)$ {\em Idem with $M=10^8$}.}
\end{center}
\end{figure}

\bigskip
\ni $\rhd$  {\sc Up \& out Call option:} Let $0\le K\le L$. The Up-and-Out Call option with strike $K$ and barrier $L$ is defined by
its payoff functional 
\[
F(x)= e^{-rT}\left(x(T)-K \right)_+\mbox{\bf 1}_{\{\max_{s\in [0,T]}x(s)\le L\}}, \quad x\!\in {\cal C}([0,T],\R).
\]

It is again classical background, that in a $B$-$S$ model
\begin{eqnarray*}
{\rm Call}^{U\&O}(X_0,r,\sigma, T) \!\!&\!\!=\!\!&\!\! {\rm Call}^{BS}\!(X_0,K,r,\sigma,T)- {\rm
Call}^{BS}\!(X_0,L,r,\sigma,T) -e^{-rT}(L\!-\!K)\Phi(d^-(L))\\
\!\!&\!\!\!\!&\!\!-\left(\!\frac{L}{X_0}\!\right)^{1+\mu}\hskip -0.55 cm \left({\rm
Call}^{BS}\!(X_0,K',r,\sigma,T)\!-\! {\rm Call}^{BS}\!(X_0,L',r,\sigma,T)
\!-\!e^{-rT}(L'\!-\!K')\Phi(d^-(L'))\right)
\end{eqnarray*}
with \[
K'= K\left(\frac{X_0}{L}\right)^2,\; L'= L\left(\frac{X_0}{L}\right)^2,\; d^-(L) =
\frac{\log(X_0/L)+(r-\frac{\sigma^2}{2})T}{\sigma\sqrt{T}}
\; \mbox{ and }\; \Phi(x):=\int_{-\infty}^x
\!\!e^{-\frac{\xi^2}{2}}\frac{d\xi}{\sqrt{2\pi}}
\]
and $\mu =\displaystyle \frac{2r}{\sigma^2}$.

We took again the same values for the $B$-$S$ parameters as for the vanilla call. We set the barrier value  at
$L=300$.   For this set of parameters $C^{UO}_0=8.54$. We tested the same three schemes. The numerical results are depicted in
Figure~5.

The conclusion (see Figure~5(a) (left))  is that, at this very high level of volatility, when $M=10^6$ 
(which is a standard size given the  high volatility setting) the (quasi-)consistent
$3$-step
$R$-$R$ extrapolation with Brownian bridge clearly outperforms the continuous Euler scheme (Brownian bridge)  of equivalent complexity
while the $3$-step $R$-$R$ extrapolation based on the stepwise constant Euler schemes with consistent Brownian increments is not
competitive at all: it suffers from both a too high variance (see Figure~5(a) (right)) for the considered sizes of the Monte Carlo
simulation and from its too slow rate of convergence in time.  

When $M=10^8$  (see Figure~5(b) (left)), one verifies again that  the time discretization error of the $3$-step $R$-$R$
extrapolation almost vanishes like for the partial lookback option. This no longer the case with the $3$-step $R$-$R$ extrapolation of
stepwise constant Euler schemes. It seems clear that the discretization time error is more prominent for the barrier  option: thus
with $n=10$, the  relative error is $\frac{9.09-8.54}{8.54}\approx 6.5 \%$ by this first $R$-$R$ extrapolation  whereas, the $3$-step
$R$-$R$ method based on the quasi-consistent  ``bridged" method yields a an approximate premium of $8.58$ corresponding to a relative
error of 
$\frac{8.58-8.54}{8.54}\approx 0.4\%$. These specific results (obtained without any control variate) are representative of the global
behaviour of the methods as emphasized by Figure~5(b)(left).

%$9.0844$.
\begin{figure}
\begin{center}
\begin{tabular}{lcc}
$(a)$&&\\ $\includegraphics[height =5cm,width=7.5cm]{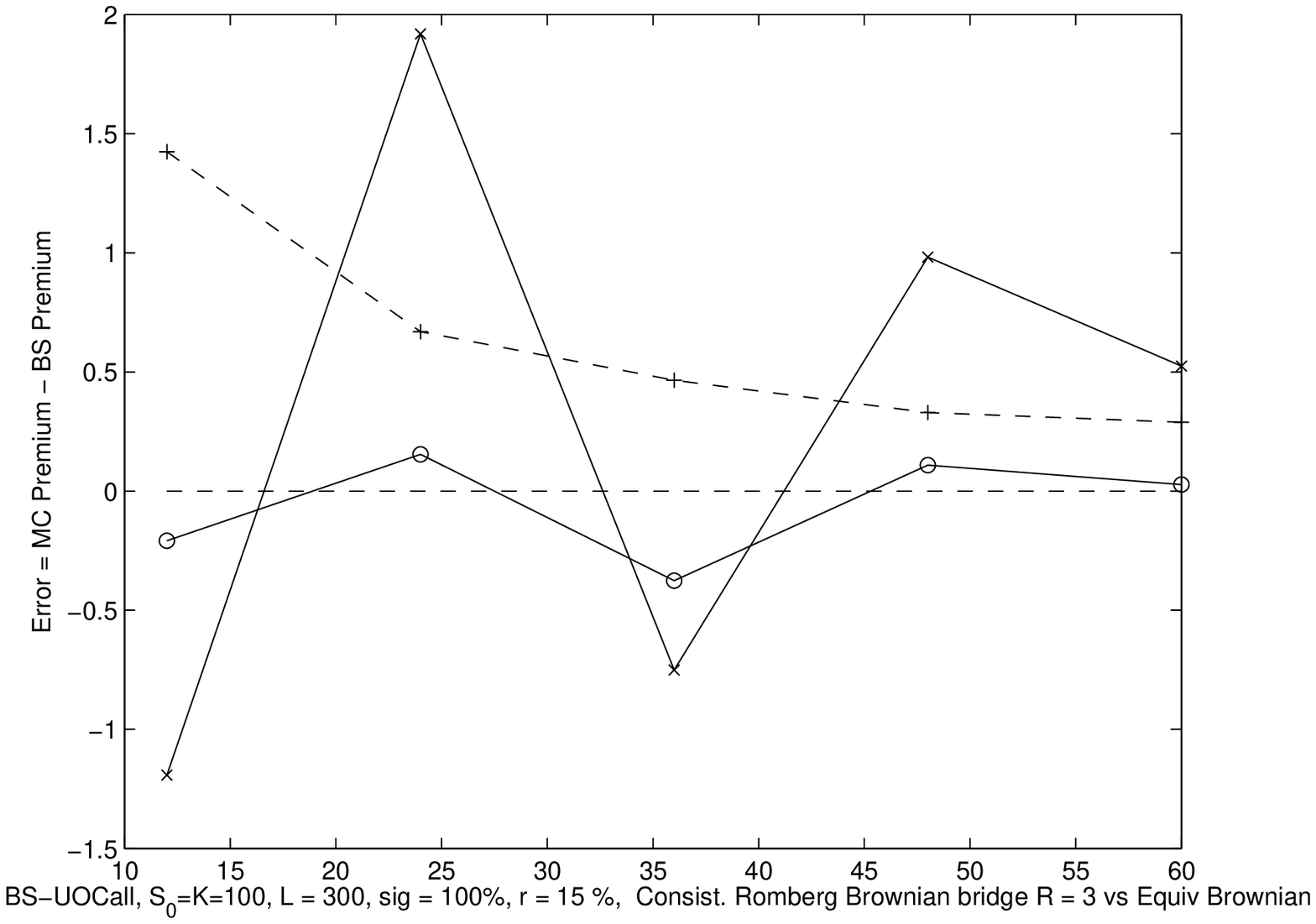}$&
$\includegraphics[height =5cm,width=7.5cm]{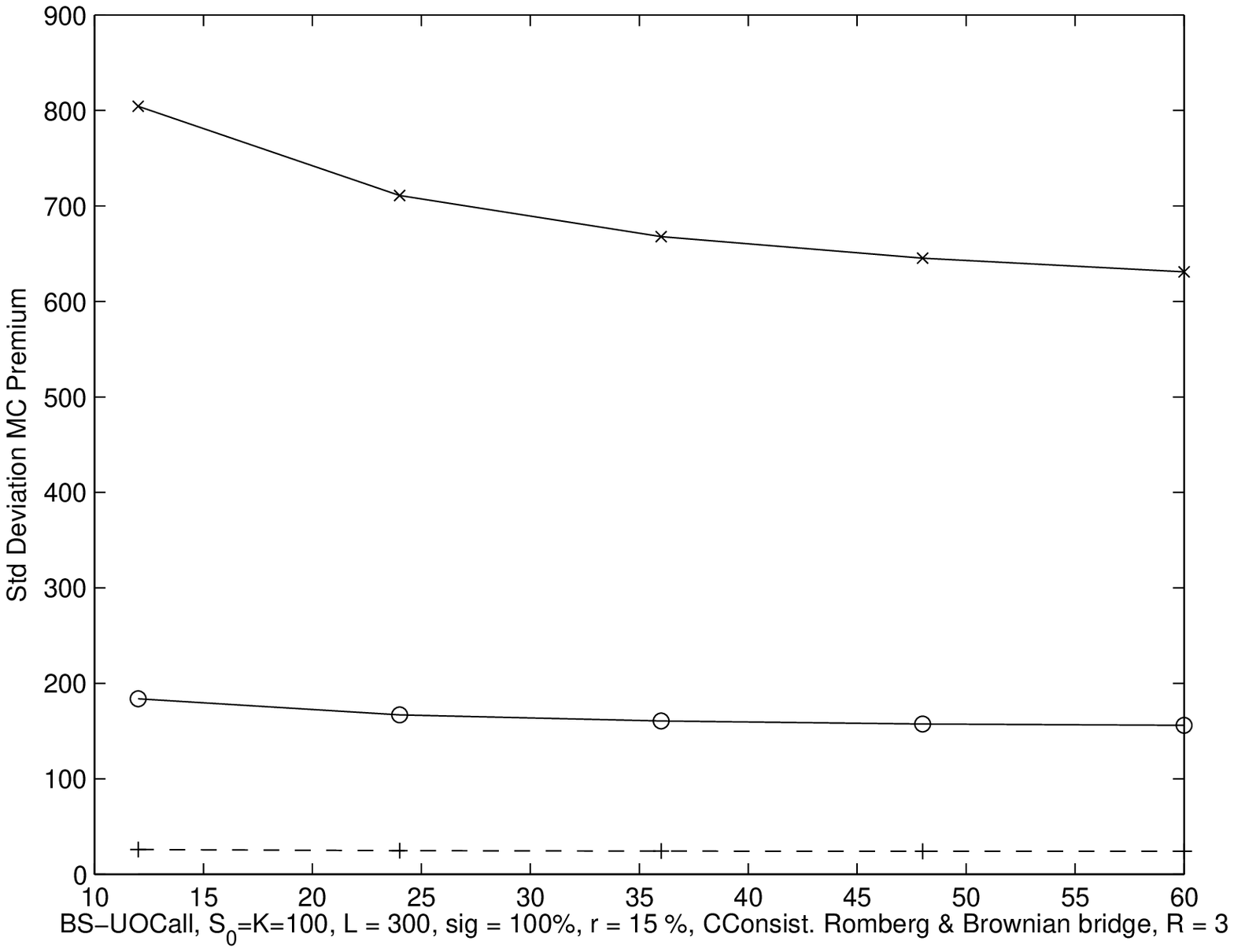}$\\
$(b)$&&\\ $\includegraphics[height =5cm,width=7.5cm]{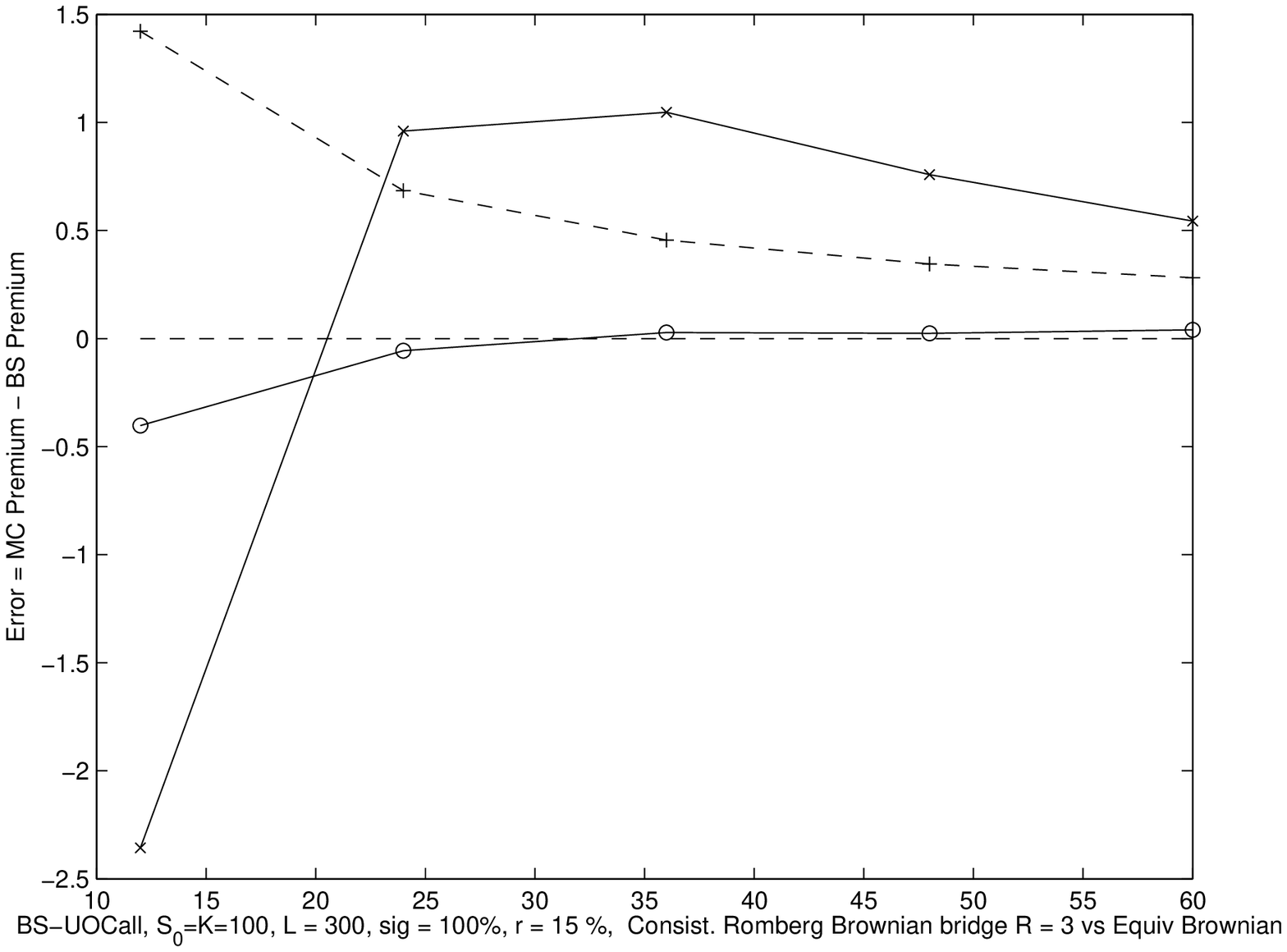}$&
$\includegraphics[height =5cm,width=7.5cm]{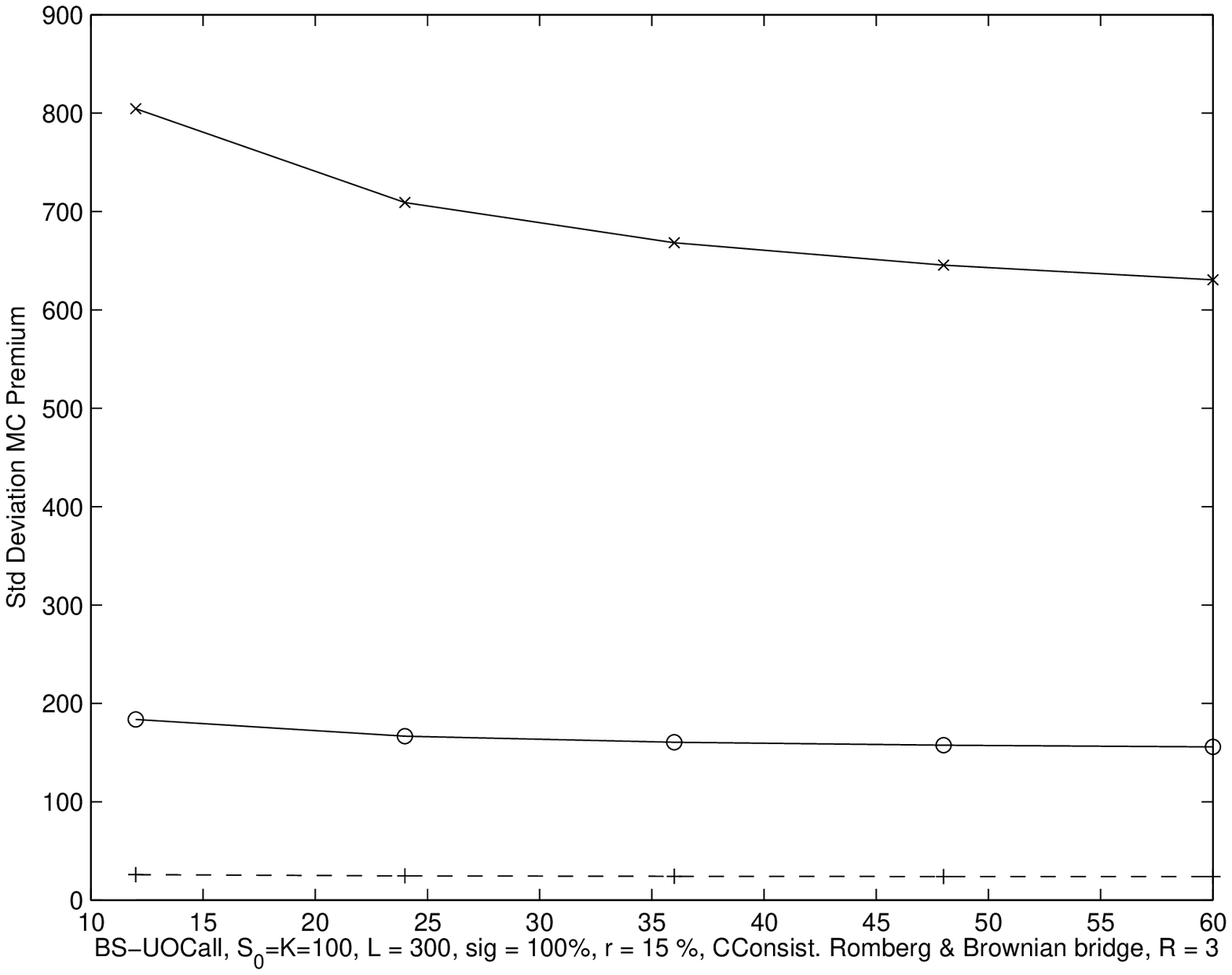}$
%&$$\includegraphics[height =4 cm,width=4cm]{Fig3(erg).eps}$$
 \end{tabular}
\caption{\small{\sc   B-S Euro up-\&-out Call option.} $(a)$ $M\!=\!10^6$. {\em $R$-$R$ extrapolation  ($R\!=\!3$) of the Euler scheme
with Brownian bridge: $\!o\!-\!\!\!-\!o\!-\!\!\!-\!o\!$.   Consistent  $R$-$R$ 
extrapolation  ($R\!=\!3$): ---$\!\times\!$---$\!\times\!$---$\!\times\!$---. Euler scheme with Brownian bridge and equivalent
complexity: $+- -+-  -+$. $X_0\!=\!K\!=\!100$, $L=300$, $\sigma\!=\!100 \%$, $r\!=\!15
\%$. Abscissas: $6n$, $n=2,4,6,8,10$. 
 {\rm Left:} Premia. {\rm Right:} Standard Deviations.} $(b)$ {\em Idem with $M=10^8$}.}
\end{center}
\end{figure}

\section{Conclusion} The multi-step $R$-$R$ extrapolation method with consistent Brownian increments provides an efficient method to
evaluate expectations of functionals of diffusions having a very high diffusion coefficient using a simple Monte Carlo simulation
based on stepwise constant or continuous Euler schemes of reasonable size (in terms of discretization). This is made possible by the control of both
  variance and    complexity.  

However the asymptotic variance control may have not produced its effect when the time discretization parameter $n$ is too small.
Then it could be useful to explore some on-line variance reduction method: the idea would be to use some stochastic approximation  methods as
introduced in~\cite{AR} to specify directly the optimal variance structure of the Euler schemes involved in the extrapolation
 for a given value of $n$.

\bigskip

\small
\noindent {\sc Acknowledgment:} We thank Julien Guyon for fruitful comments on  a preliminary version and \'Eric Sa\"{\i}as for his help on
Number Theory results.

\normalsize 

\bigskip
\centerline{\Large \bf Annex}

\medskip

We will show that in some situations the equality $X^{(1)}_{_T}= X^{(2)}_{_T}$ $a.s.$ may imply that $W^{(1)}=W^{(2)}$ $\P$-$a.s.$ as processes
defined on
$[0,T]$. Assume that $(W^{(1)}, W^{(2)})$ is a   $\R^{2q}$-dimensional Brownian motion with respect to its (augmented) natural
filtration  denoted $({\cal F}^{1,2}_t)_{t\in [0,T]}$, both marginals $W^{(1)}$ and $W^{(2)}$ being standard $d$-dimensional Brownian motions.
Let $R_{_W}:=\left[\E (W^{(1)}_1W^{(2)*}_1)\right]\!\in{\cal M}(q\times q)$ denote the correlation matrix of 
$W_1^{(1)}$ and $W_1^{(2)}$.

\begin{Proposition}  
Let $F_{b,\sigma}: [0,T]\times \R^d\times \R^d \to \R$ denote the function defined by 
\[
F_{b,\sigma}(t,x_1,x_2) = (x_1-x_2)^*(b(t,x_1)-b(t,x_2))+\frac 12 \sum_{i=1}^d \left(|\sigma_{i.}(t,x_1)|-|\sigma_{i.}(t,x_2)|\right)^2.
\]
Assume that, for every $t\!\in [0,T]$,  
\begin{equation}\label{Fbsigma}
\forall\, x_1,\, x_2 \!\in \R^d,  \qquad x_1\neq x_2  \;\Longrightarrow\;  F_{b,\sigma}(t,x_1,x_2)>0. 
\end{equation}
Furthermore, assume 
\begin{equation}\label{noisy}
\exists i\!\in \{1,\ldots,d\},\quad \forall\, j\!\in \{1,\ldots,q\},\qquad \int_0^T \E\left(\sigma^2_{ij}(t,
X_t)\right)dt >0.
\end{equation}

\smallskip
If $X^{(1)}_{_T}= X^{(2)}_{_T}$ $\P$-$a.s.$, then  $W^{(1)}= W^{(2)}$ (and $X^{(1)}=X^{(2)}$) $\P\mbox{-}a.s.$.  
\end{Proposition}

\noindent {\bf Remark.} Assumption~(\ref{Fbsigma}) is always satisfied if $b$ is increasing $i.e.$ 
\[
\forall\, x_1,\, x_2 \!\in \R^d,  \qquad x_1\neq x_2\Longrightarrow  (x_1-x_2)^*(b(t,x_1)-b(t,x_2))>0.
\]

\bigskip
\noindent {\bf Proof.}  $(a)$ It follows from It\^o's
formula that 
\begin{eqnarray*}
 |X^{(1)}_{_T}-X^{(2)}_{_T}|^2&=& \int_0^T\Phi(t,X^{(1)}_{t},X^{(2)}_{t})dt +M_{_T}
\end{eqnarray*}
where 
\begin{eqnarray*}
\Phi(t,x_1,x_2)&=& 2(x_1-x_2)^*\,(b(t,x_1)-b(t,x_2))+{\rm
Tr}(\sigma\sigma^*(t,x_1) )+{\rm Tr}(\sigma\sigma^*(t,x_2))\\
&& -2 
\sum_{i=1}^d\sigma_{i.}(t,x_1)^*\,R_{_W}\sigma_{i.}(t,x_2)
\end{eqnarray*}
and $M_t=2\int_0^t (X^{(1)}_{s}-X^{(2)}_{s})(\sigma(t,X^{(1)}_{s})dW^{(1)}_s-\sigma(t,X^{(2)}_{s})dW^{(2)}_s)$ is an ${\cal F}_t$-local
martingale null at zero. In fact it is a true martingale since all the coefficients have linear growth and $\sup_{t\in [0,T]}|X^{(1)}_t|$ lies in
every
$L^p(\P)$, $p>0$. One checks from the definition of $R_{_W}$ that, for
every
$u,\,v\!\in \R^q$, $u^*\,R_{_W}v
\le |u|\,|v|$. 
Consequently, for every $t\!\in [0,T]$, every  $x_1,\, x_2\!\in \R^d$,  
\[
\Phi(t,x_1,x_2) \ge 2F_{b,\sigma}(t,x_1,x_2)\ge 0.
\]
If $X^{(1)}_{_T}= X^{(2)}_{_T}$ $\P$-$a.s.$, then 
$$
M_{_T}=-\int_0^T \Phi(t,X^{(1)}_{t},X^{(2)}_{t})dt\le 0\qquad \mbox{$\P$-$a.s.$}.
$$ 
Hence  $M_{_T}= 0$ $\P$-$a.s.$ since $\E M_{_T}=0$. In turn this implies that 
$\int_0^T\Phi(t,X^{(1)}_{t},X^{(2)}_{t})dt=0$ $\P$-$a.s.$. The  continuous  function $\Phi$ being non negative  and $(X^{(1)}_{t},X^{(2)}_{t})$
being pathwise continuous, 
\[
\P\mbox{-}a.s. \qquad \forall\, t\in [0, T],\qquad \Phi(t,X^{(1)}_{t},X^{(2)}_{t})= F_{b,\sigma}(t,X^{(1)}_{t},X^{(2)}_{t})=0.
\]
Consequently, (\ref{Fbsigma}) 
\[
\P\mbox{-}a.s. \qquad\forall\, t\!\in [0,T], \qquad X^{(1)}_{t}=X^{(2)}_{t}
\]

Elementary computations show that 
\[
\Phi(t,\xi,\xi) =  2 \sum_{i=1}^d  \sigma_{i.}(t,\xi)^*\,(I_q-R_{_W} )\sigma_{i.}(t,\xi).
\]
The symmetric matrix $I_q-R_{_W}$ being nonnegative, 
\[
\forall\, i\!\in \{1,\ldots,d\}, \quad \sigma_{i.}(t, X^{(1)}_{t})^*(I_q-R_{_W} )\sigma_{i.}(t, X^{(1)}_{t})=0\quad i.e.\quad\sigma_{i.}(t,
X^{(1)}_{t})\!\in {\rm Ker}(I_q-R_{_W}).
\]
If $I_q\neq R_{_W} $, this (nonnegative symmetric) matrix has at least one positive eigenvalue $\lambda>0$. Let $u\!\in \R^q\setminus\{0\}$ be an
eigenvector  associated to
$\lambda$. Then $u\perp  {\rm Ker}(I_q-R_{_W})$   so that $u^*\,\sigma_{i.}(t, X^{(1)}_{t})=0$ $\P$-$a.s.$. This cannot be satisfied
by an index  $i$ satisfying~(\ref{noisy}). Hence $R_{_W}=I_q$ which in turn implies that $W^{(1)}=W^{(2)}$ since $(W^{(1)}, W^{(2)})$ is a Gaussian
centered process.$\qquad_\diamondsuit$

\end{document}